\documentclass{article}
\usepackage{graphicx}%
\usepackage{multirow}%
\usepackage{amsmath,amssymb,amsfonts}%
\usepackage{amsthm}%
\usepackage{mathrsfs}%
\usepackage[title]{appendix}%
\usepackage{xcolor}%
\usepackage[export]{adjustbox}
\usepackage{bm}
\usepackage{footnotehyper}
\usepackage{amssymb,amsmath,amsthm}
\usepackage{amscd}
\usepackage{hyperref}
\usepackage{lipsum}
\usepackage{epsfig}
\usepackage{amsfonts}
\usepackage{amssymb}
\usepackage{amsmath,enumerate}
\usepackage{commath}
\usepackage{euscript}
\usepackage{enumitem}
\usepackage{amscd}
\usepackage{xcolor}
\usepackage[numbers,sort&compress]{natbib}
\usepackage{lscape}
\usepackage{textcomp}%
\usepackage{manyfoot}%
\usepackage{booktabs}%
\usepackage{algorithm}%
\usepackage{algorithmicx}%
\usepackage{algpseudocode}%
\usepackage{listings}%
\usepackage{multirow}
\usepackage{amsmath, amsthm, amssymb}
\usepackage{graphicx}
\usepackage{algorithm,setspace}
\usepackage{rotating}
\usepackage{enumitem}
\usepackage{algpseudocode}
\newtheorem{corollary}{Corollary}[section]  
\usepackage{commath}
\usepackage{booktabs} 

\newtheorem{lemma}{Lemma}  
\newtheorem{assumption}{A}
\usepackage{hyperref}
\hypersetup{
    colorlinks=true,
    linkcolor=blue,
    filecolor=blue,
    urlcolor=blue,
}
\theoremstyle{thmstyleone}%
\newtheorem{theorem}{Theorem}%  meant for continuous numbers
\newtheorem{proposition}[theorem]{Proposition}% 
\theoremstyle{thmstyletwo}%
\newtheorem{example}{Example}%
\newtheorem{remark}{Remark}%

\theoremstyle{thmstylethree}%
\newtheorem{definition}{Definition}%

\begin{document}
\title{Iterative Methods for Projected Solutions of Quasi-equilibrium Problems}
\author{Didier Aussel\footnotemark[1]~\footnotemark[2],~ 
    Jauny\footnotemark[3],~  
    Asrifa Sultana\footnotemark[3],~  
    Shivani Valecha\footnotemark[3]}

    \renewcommand{\thefootnote}{\fnsymbol{footnote}}
\maketitle

\footnotetext[1]{Corresponding author. Email: \texttt{aussel@univ-perp.fr}}
\footnotetext[2]{CNRS PROMES, UPR 8521, Université de Perpignan Via Domitia, Perpignan, France.}
\footnotetext[3]{Department of Mathematics, Indian Institute of Technology Bhilai, 491002, India.}

%\footnotetext[1]{ Corresponding author. e-mail- {\tt asrifa@iitbhilai.ac.in}}
%\noindent
%\footnotetext[2]{Department of Mathematics, Indian Institute of Technology Bhilai, Raipur - 492015, India.
%}
%  \author*[1]{\fnm{Didier} \sur{Aussel}}\email{aussel@univ-perp.fr}
 
%  \author[2]{\fnm{} \sur{Jauny}}\email{jaunyp@iitbhilai.ac.in}
% % \equalcont{These authors contributed equally to this work.}
% \author[2]{\fnm{Asrifa} \sur{Sultana}}\email{asrifa@iitbhilai.ac.in}
% \author[2]{\fnm{Shivani} \sur{Valecha}}\email{shivaniv@iitbhilai.ac.in}
% \affil*[1]{\orgdiv{CNRS PROMES, UPR 8521}, \orgname{Université de Perpignan Via Domitia}, \orgaddress{ \street{Technosud}, \city{Perpignan}, \country{France}}}

% \affil[2]{\orgdiv{Department of Mathematics}, \orgname{Indian Institute of Technology}, \orgaddress{ \city{Bhilai}, \postcode{491002}, \state{Chhattisgarh}, \country{India}}}

% \affil*[1]{\orgdiv{Department}, \orgname{Organization}, \orgaddress{\street{Street}, \city{City}, \postcode{100190}, \state{State}, \country{Country}}}

% \affil[2]{\orgdiv{Department}, \orgname{Organization}, \orgaddress{\street{Street}, \city{City}, \postcode{10587}, \state{State}, \country{Country}}}

% \affil[3]{\orgdiv{Department}, \orgname{Organization}, \orgaddress{\street{Street}, \city{City}, \postcode{610101}, \state{State}, \country{Country}}}

\begin{abstract}
Aussel et al. (J Optim Theory Appl 170:818–837, 2016) introduced the concept of projected solutions for the quasi-variational inequalities with a non-self constraint map, that is, the case where the constraint map may take values
outside the feasible set. This paper presents two iterative methods to determine the projected solutions for quasi-equilibrium problems (QEP). We prove the convergence of the sequence generated by iterative methods to a projected solution of the QEP under some suitable assumptions. Some encouraging numerical experiments are presented to show the performance of the proposed methods. As an application, we apply the proposed algorithms to solve an electricity market model.
\end{abstract}

%\keywords{keyword1, Keyword2, Keyword3, Keyword4}

% In the proposed methods, a projected solution for QEP is computed by solving a regularized equilibrium problem at each iteration and then successively projecting it on the feasible set.
\section{Introduction}

\noindent
 The concept of equilibrium problems, introduced by Blum and Oettli \cite{blum}, has gained a lot of attention because it includes many important problems such as variational inequalities, complementarity problems, optimization problems, fixed point problems, and saddle point problems. For a closed and convex subset \(\mathcal{D} \subset \mathbb{R}^n\) and  bifunction \(f : \mathcal{D} \times \mathcal{D} \rightarrow \mathbb{R}\), the aim of equilibrium problem is to determine \(\bar{y} \in \mathcal{D}\) such that 
 \begin{equation}\label{EP}
    f(\bar{y}, y) \geq 0, \quad \forall\,y \in \mathcal{D}.
 \end{equation}
  The equilibrium problem \eqref{EP} and its set of solutions are indicated by \(\text{EP}(f, \mathcal{D})\) and $\text{S}_{\text{EP}}(f, \mathcal{D})$, respectively. Further, if EP($f,\mathcal{D}$) admits a unique solution, then we denote this solution by $\bar y(\mathcal{D})$. These problems have been studied in recent years because of their wide applications (see e.g., \cite{quoc,book_eq_prob,cotrina_eq_prob,bigi}). %The dual problem related to \(\text{EP}(f, \mathcal{C})\) is to find \(y^* \in \mathcal{C}\) such that
%\begin{equation}\label{dual_EP}
%    f(x, y^*) \leq 0, \quad \forall x \in \mathcal{C}
%\end{equation}
% and its solution set is denoted by \(\text{S}^{d}_{\text{EP}}(f, \mathcal{C})\).

The quasi-equilibrium problem (QEP) is an extension of the equilibrium problem \eqref{EP} in which the constraint set depends on the current point. More precisely, the aim of the QEP is to find \(\bar{x} \in \mathcal{K}(\bar{x})\) such that
\begin{equation}\label{QEP}
   f(\bar{x}, y) \geq 0, \quad\forall\, y \in \mathcal{K}(\bar{x}), 
\end{equation}
where \(\mathcal{K}:\mathcal{C}\rightrightarrows\mathbb{R}^n\) is the constraint map for a given set \(\mathcal{C}\subset \mathbb{R}^{n}\). The QEP \eqref{QEP} and its solution set will be denoted by QEP($f, \mathcal{K}$) and \(S_{\text{QEP}}(f, \mathcal{K})\), respectively. A variety of mathematical problems can be reformulated as quasi-equilibrium problems, such as quasi-optimization problems, generalized Nash equilibrium problems, and quasi-variational inequalities (see, for instance, \cite{Aussel_book,Cotrina_note_QEP} and references therein). In particular, for a given mapping \(F: \mathcal{C} \rightrightarrows \mathbb{R}^{n}\) with compact values if 
\begin{equation}\label{Inro_QVI_f}
f(x, y) = \sup_{{x^*}\in F(x)}\langle {x^*}, y - x \rangle
\end{equation}
then QEP($f, \mathcal{K}$) reduces to quasi variational inequality QVI($F, \mathcal{K}$). The existence results for QEP($f, \mathcal{K}$) and QVI($F, \mathcal{K}$) are studied in the literature under the assumption that the constraint map $\mathcal{K}$ is a self-map, i.e., $\mathcal{K}(\mathcal{C})\subseteq \mathcal{C}$ \cite{Cotrina_note_QEP, Aussel_Asrifa}.

%If \(f(x, y) = \sup_{x^*\in F(x)}\langle x^*, y - x \rangle\), for a given mapping \(F: \mathcal{C} \rightrightarrows \mathbb{R}^{n}\) with compact values, EP reduces to the variational inequality problem (VIP) and QEP reduces to the quasi-variational inequality problem (QVIP).
In \cite{aussel2016existence}, Aussel et al. studied the problem QVI($F, \mathcal{K}$) where the constraint map is not a self-map, i.e., $\mathcal{K}(\mathcal{C})\nsubseteq \mathcal{C} $. Since the existence of a point $x\in \mathcal{K}(x)$ satisfying QVI($F, \mathcal{K}$) cannot be guaranteed in this case, the authors in \cite{aussel2016existence} initiated the concept of projected solutions for these problems motivated by an illustration of an electricity market model. Later, Cotrina and Z{\'u}{\~n}iga \cite{cotrina2019quasi} adapted this notion of projected solutions for the quasi-equilibrium problems with non-self constraint maps and extended the existence results in \cite{aussel2016existence} to the case of QEP. It is worth mentioning that the existence of projected solutions for various problems like QEP, QVI and generalized Nash games is well studied in the recent past \cite{sara,milasi_proj,cotrinaGNEP}. 

On the other hand, numerous algorithms have been formulated in the last few decades to find the solution of QVI and QEP when the constraint map $\mathcal{K}$ is a self-map (e.g., \cite{scrimalli,Bigi G,QEP_proximal}). Nesterov-Scrimali \cite{scrimalli} formulated an algorithm for quasi-variational inequalities in which certain strongly convex optimization problems are solved in each iteration. The convergence of this algorithm to a classical solution of QVI is shown by assuming that the constraint map $\mathcal{K}:\mathbb{R}^n\rightrightarrows \mathbb{R}^n$ satisfies a contraction-type property (see \cite[Theorem 4.4]{scrimalli}). %devised an algorithm for variational inequality problems in which the dual gap function is used and strongly convex optimization problems are solved in each iteration. Based on this algorithm for VI problems, the authors in \cite{scrimalli} further . 
Santos-Souza \cite{QEP_proximal} studied the proximal point algorithm for QEP in which an equilibrium problem is solved in each iteration by regularizing the given bifunction. The convergence of this algorithm to a classical solution of QEP is shown under the assumptions that $\mathcal{K}:\mathcal{C}\rightrightarrows \mathcal{C}$ is continuous and a special solution set of QEP is non-empty (see \cite[Theorem 3.6]{QEP_proximal}). Recently, Bianchi et. al. \cite{bianchi2025projected} formulated an algorithm for finding the projected solution of QEP where the constraint map need not be a self-map. The convergence of this algorithm to a projected solution is shown by assuming that the map $\mathcal{K}:\mathcal{C}\rightrightarrows \mathbb{R}^n$ is continuous and the generated sequence is asymptotically regular \cite[Theorem 3.1]{bianchi2025projected}. %The convergence of these algorithms to a solution is strongly dependent on the condition $\mathcal{K}(\mathcal{C})\subseteq \mathcal{C}$. %When the constraint map $\mathcal{K}$ is a self-map, numerous algorithms based on the extragradient-type methods, Newton-type method, augmented Lagrangian method and gap functions technique are available in the literature to solve QEP$(f, \mathcal{K})$, for example, see \cite{Strodiot,Van NTT,Santos PJS,Bueno LF, Bigi G} and references therein. For solving QEP$(f,\mathcal{K})$, Newton-type and augmented Lagrangian approaches are very efficient. However, they need the problem to be differentiable. Furthermore, extragradient-type methods, Newton-type methods, and augmented Lagrangian methods require a good choice of initial point. Therefore, the aim of this paper is to generalize those iterative methods to find the projected solutions of QEP$(f,\mathcal{K})$ that do not have these drawbacks. So far, an algorithm to find the projected solution of QEP can only be found in \cite{bianchi2025projected}.

In this paper, we consider quasi-equilibrium problem in which the constraint map $\mathcal{K}$ need not be self, that is, $\mathcal{K}(\mathcal{C})\nsubseteq \mathcal{C}$. We propose two iterative methods to find the projected solutions for the considered QEP. In our first algorithm, we solve certain strongly convex optimization problems in each iteration, motivated by \cite{scrimalli}. %is an extension of the method proposed in \cite{scrimalli} for the classical solution of QVI. In this iterative algorithm, 
 We show the convergence of this algorithm to a projected solution under the strong monotonicity and Lipschitz-type continuity assumption on the bifunction $f$. In our second proposed algorithm, which is inspired by the proximal point method \cite{QEP_proximal}, we need to solve a strongly monotone EP by regularizing the bifunction in each iteration. %and  second algorithm for finding the projected solution of QEP is an extension of the proximal point algorithm studied in . 
It can be noted that we have shown the convergence of the proximal point method for QEP without requiring the strong monotonicity or Lipschitz-type continuity conditions on the bifunction $f$. However, we require the nonemptiness of the special solution set for the considered QEP in this algorithm. Our algorithms are applicable to generalized Nash equilibrium problems where the constraint map $\mathcal{K}$, formed by the product of feasible strategy maps of players, need not be a self-map. We apply our algorithms to find a projected solution for the deregulated electricity market model studied in \cite{aussel2016existence}. We also illustrate our methods with some encouraging numerical experiments.

%The convergence of the algorithm proposed in \cite{bianchi2025projected} is proved by assuming that the sequence produced is asymptotically regular, which is not easy to verify in practice.

%So far,  Despite its complexity, depending strongly on the method used to solve the subproblems, the simplicity of implementation and the fact that each iterate of the method is a solution of an equilibrium problem make the proximal point method quite attractive for solving quasi-equilibrium problems. 
 %Furthermore, it can be performed for any initial point, and there are a larger number of algorithms available for EP than for QEP to solve the subproblems in each iteration of the method

 The remaining paper is structured as follows. After presenting some basic concepts and notation in quasi-equilibrium problems in Section \ref{prelim}, we propose two iterative algorithms and establish their respective convergence properties in Section \ref{iter_algos}. In Section \ref{num_exp}, some numerical results are reported to show the performance of the proposed methods. In Section \ref{sec_GNEP}, we provide the sufficient conditions to find the projected solutions for the generalized Nash equilibrium problems by using the proposed algorithms. Finally, we apply the proposed algorithms to calculate the projected solutions for an electricity market model in Section \ref{EMM_NUM}.

\section{Basic Concepts}\label{prelim}
Let $\overline{co}(P)$ be the closure of the convex hull of $P\subseteq \mathbb{R}^n$. For the constraint map $\mathcal{K}:\mathcal{C}\rightrightarrows \mathbb{R}^n$ considered in QEP($f,\mathcal{K}$) \eqref{QEP}, we denote $\mathcal{A}=\overline{co}(\mathcal{K}(\mathcal{C}))$. For any $r\in \mathbb{R},$ the expression $\lceil r \rceil$ represents the least integer greater than or equal to $r$. %$\mathcal{C} \subseteq \mathbb{R}^{n}$ is a nonempty closed convex set and $\mathcal{A}=\overline{co}(\mathcal{K}(\mathcal{C})).$
We recall the notion of projected solutions for a quasi-equilibrium problem in which the constraint map is not necessarily a self-map.

\begin{definition}\cite{cotrina2019quasi}
 Let $\mathcal{K}: \mathcal{C} \rightrightarrows \mathbb{R}^{n}$ be a set-valued map and $f: \mathcal{A} \times \mathcal{A} \to \mathbb{R}$ be a bifunction. A point $\bar{x}  \in \mathcal{C}$ is said to be a \textit{projected solution} of the QEP$(f, \mathcal{K})$ if there exists $\bar{y}\in \mathcal{K}(\bar{x})$ such that
\begin{itemize}
\item[(a)] $\norm{\bar{x}-\bar{y}}=\inf_{y\in \mathcal{C}}  \norm{\bar{x}-{y}}$;
\item[(b)]\label{def_b} $f(\bar y,y)\geq 0,~\text{for all}~y\in  \mathcal{K}(\bar{x}).$
	\end{itemize}  
\end{definition}
\noindent The set of projected solutions will be denoted by $S^{P}_{QEP}(f, \mathcal{K}).$ Clearly, if $\bar{y}\in \mathcal{K}(\bar{x})\cap\mathcal{C}$ satisfies the above inequality (b) in Definition \ref{def_b}, then it becomes the classical solution of the quasi-equilibrium problem. Hence, the necessary condition for the existence of a classical solution of QEP($f,\mathcal{K}$) \eqref{QEP} is $\mathcal{K}(\mathcal{C})\cap \mathcal{C} \neq \emptyset.$ %This can only happen when $\mathcal{K}(\mathcal{C})\cap \mathcal{C} \neq \emptyset.$ %$\mathcal{K}$ is a self-map, i.e., $\mathcal{K}(\mathcal{C})\subseteq \mathcal{C}$. When this is not the case, then either $\mathcal{K}(\mathcal{C})\cap \mathcal{C} = \emptyset$ or  $\mathcal{K}(\mathcal{C})\cap \mathcal{C} \neq \emptyset.$
%Throughout this paper, we denote \(\| x \| = \sqrt{x^\top x}\), for all \(x \in \mathbb{R}^n\). %Let \(\mathcal{C} \subset \mathbb{R}^n\) be a closed and convex set, and a bifunction \(f : \mathcal{C} \times \mathcal{C} \rightarrow \mathbb{R}\) satisfying the following conditions:

%strongly monotone if $f(x,y)+f(y,x) \le -\mu \|x-y\|^2$ for all $x,y \in \mathcal{A}$ and $\mu>0.$ 

%Now,  The following examples show the possible relationships between classical and projected solutions in these situations.
For $X\subseteq \mathbb{R}^n$, define a projection map $P_X:\mathbb{R}^n\rightrightarrows X$ as,
\begin{equation*}\label{Pr}
	P_X(y)=\{x\in X\,|\,\norm{y-x}=\inf_{w\in X}\norm{y-w}\}.
\end{equation*}
If $X$ is a non-empty closed convex set, then the projection map $P_X$ defined as above is a single-valued, continuous map and satisfies the non-expansivity property \cite{stampbook}, that is,
\begin{equation*}
    \norm{P_X(x)-P_X(y)}\leq \norm{x-y},~~\forall\,x,y\in \mathbb{R}^n.
\end{equation*}
Further, in this case, $P_X(y)=\{x\}$ if and only if,
	\begin{equation}\label{projection_inequality}
		\langle x-y, \eta-x\rangle \geq 0,\quad\forall\,\eta\in X.
	\end{equation}
%  Let us recall the following result related to the projection map, which will be used to prove upcoming equivalence and existence results on projected solutions.
%\begin{lemma}\label{projection}\cite[Lemma 2.1 and Theorem 2.3]{stampbook}
%	Let $X\subset\mathbb{R}^n$ be non-empty closed convex set and the map $P_X$ be defined as (\ref{Pr}). Then, for each $y\in \mathbb{R}^n$ there exists a unique $x\in X$ such that $P_X(y)=\{x\}$. Further, $P_X(y)=\{x\}$ if and only if,
	%\begin{equation}\label{projection_inequality}
	%	\langle x-y, \eta-x\rangle \geq 0,~\text{for all}~\eta\in X.
	%\end{equation}
%\end{lemma}

Let us recall some concepts related to the monotonicity of bifunctions. For any set $\mathcal{D}\subseteq\mathbb{R}^n$, a bifunction $f:\mathcal{D}\times \mathcal{D} \to \mathbb{R}$ is said to be 
\begin{itemize}
    \item[-] \textit{monotone} if $f(x,y) + f(y,x) \le 0$, for all $x,y \in \mathcal{D}$;
    \item[-] \textit{strictly monotone} if $f(x,y) + f(y,x) < 0$, for all $x,y \in \mathcal{D}$;
    \item[-]\textit{strongly monotone} with parameter $\mu>0$ if the following inequality holds:
    \begin{equation*}\label{strongly_mon}
    f(x,y) +f(y,x) \le -\mu \|x-y\|^2,\quad\forall\, x,y\in \mathcal{D}.    
    \end{equation*}        
\end{itemize}
It is easy to observe that strong monotonicity $\Rightarrow$ strict monotonicity $\Rightarrow$ monotonicity.
%will be helpful in a further analysis of the proposed algorithms.

% \begin{definition}\cite{quoc}
%     monotone if the following condition holds:
%     \begin{equation}\label{monotone_bifun}
%        f(x,y) + f(y,x) \le 0~~\forall\,x,y \in \mathcal{D}. 
%     \end{equation}
%     If the above inequality \eqref{monotone_bifun} holds strictly, then f is said to be strictly monotone. Furthermore, $f$ is said to be strongly monotone on $\mathcal{D}$ with parameter $\mu>0$ if the following inequality hold:
%     \begin{equation}\label{strongly_mon}
%     f(x,y) +f(y,x) \le -\mu \|x-y\|^2,~~~\forall x,y\in \mathcal{D}.     
%     \end{equation}        
% \end{definition}
\begin{definition}\cite{quoc}
    For any set $\mathcal{D}\subseteq\mathbb{R}^n$, a bifunction $f:\mathcal{D}\times \mathcal{D}\to \mathbb{R}$ is said to be Lipschitz-type continuous on $\mathcal{D}$ if there exists a constant $L>0$ such that
\begin{equation}\label{Lipschitz_cond}
    f(x,y) + f(y,z) \ge f(x,z) -L\|y-x\|\|z-y\|,\quad\forall\,x,y,z\in \mathcal{D}. 
\end{equation}
\end{definition} 
Note that the Lipschitz-type condition \eqref{Lipschitz_cond} implies the following Lipschitz-type condition in the sense of Mastroeni \cite{Mastroeni}:
\begin{equation}
   f(x,y) + f(y,z) \ge f(x,z) - c_1\|y-x\|^2 -c_2\|z-y\|^2,\quad\forall\,x,y,z\in \mathcal{D},  
\end{equation}
for the given constants $c_1,c_2>0$. One can prove this by taking $c_1=\frac{L}{2r}$
and $c_2=\frac{Lr}{2}$ for any $r>0$.

The subdifferential of the bifunction $f:\mathcal{A}\times \mathcal{A} \to \mathbb{R}$ with respect to its second argument evaluated at $(x,y)\in \mathcal{A}\times \mathcal{A}$ is defined as follows:
\begin{equation}\label{subdifferential}
   \partial_{2}f(x,y) = \left\{\xi \in \mathbb{R}^n|\,f(x,y')-f(x,y)\ge \langle \xi, y'-y \rangle,\quad\forall\,y'\in \mathcal{A}\right\}. 
\end{equation}

The following result shows the relation between the Lipschitz constant of the bifunction $f:\mathcal{A}\times \mathcal{A}\to \mathbb{R}$ and the Lipschitz constant of an associated subgradient map. %QEP($f,\mathcal{K}$) and the Lipschitz constant of the principal operator $F$ of quasi-variational inequality QVI($F, \mathcal{K}$).
\begin{lemma}\label{Lemma_strong_mono}
    %Suppose $f:\mathcal{A}\times \mathcal{A}\to \mathbb{R}$ is defined as (\ref{Inro_QVI_f}) for given $F:\mathcal{A}\rightrightarrows \mathbb{R}^n$. Then $f$ satisfies the condition \eqref{Lipschitz_cond} if $F$ is $L$-Lipschitz, that is, $\norm{y^*-z^*}\leq L \norm{y-z}$ for all  $y,z\in \mathcal{A}$ and for all $y^*\in F(y), z^*\in F(z)$. In particular, if $F$ is defined as $F(x)=\partial_2 f(x,\cdot)(x)$ where $\partial_2 f$ is considered as above then,
    Suppose $f:\mathcal{A}\times \mathcal{A}\to \mathbb{R}$ satisfies Assumption $\textbf{A}$ \ref{A1}. Let $F:\mathcal{A}\rightarrow \mathbb{R}^n$ be defined as $F(x)=\partial_2 f(x,\cdot)(x)$. Then,
    \begin{itemize}
        \item [(i)] $F$ is L-Lipschitz continuous, that is, there exists $L>0$ such that $\norm{x^*-y^*}\leq L\norm{x-y}$, for all $x,y\in \mathcal{A}$ with $x\neq y$ if $f$ fulfills the Lipschitz-type continuity \eqref{Lipschitz_cond} and the monotonicity conditon.
        \item [(ii)] $F$ is $\mu$-strongly monotone, i.e., $\langle x^*-y^*,x-y\rangle\geq \mu\norm{x-y}^2$, for all $x,y\in \mathcal{A}~\text{and for all}~\,x^*\in \partial_2f(x,x),y^*\in \partial_2 f(y,y)$ if $f$ is $\mu$-strongly monotone.
    \end{itemize} 
\end{lemma}
\begin{proof}
% We observe that,
% \begin{equation*}
%     -f(x,y)-f(y,z)+f(x,z)=-\sup_{x^*\in F(x)}\langle x^*,y-x\rangle-\sup_{y^*\in F(y)} \langle y^*,z-y\rangle+\sup_{x^*\in F(x)} \langle x^*,z-x\rangle.
% \end{equation*}
% Since $F$ has compact convex values, there exists $x_1^*,x_2^*$ in $F(x)$ and $y^*\in F(y)$ such that
% \begin{align*}
%     -f(x,y)-f(y,z)+f(x,z)=&-\langle x_1^*,y-x\rangle- \langle y^*,z-y\rangle+\langle x_2^*,z-x\rangle\\
%     \leq & -\langle x_2^*,y-x\rangle -\langle y^*,z-y\rangle+\langle x_2^*,z-x\rangle\\
%     = &\langle x_2^*-y^*,z-y\rangle
%     \leq L\norm{x-y}\norm{y-z}.
% \end{align*}
% Part (i) of the result follows from the fact that solving QEP($f,\mathcal{K}$) is equivalent to solving QVI($F,\mathcal{K}$) where $F(x)=\partial_2f(x,\cdot)(x)$ (see \cite[Section 2.2]{Aussel_Dutta_Pandit}).
To prove part (i), let $x,y\in \mathcal{A}$ be arbitrary with $x\neq y$. We observe from \eqref{subdifferential} that $f(x,y)+f(y,x)\geq \langle x^*,y-x\rangle+\langle y^*,x-y\rangle$. Hence, by monotonicity of $f$, it holds that
\begin{equation*}
   -f(x,y)-f(y,x)\geq \langle x^*-y^*,y-x\rangle.
\end{equation*}
By taking $z=x$ in \eqref{Lipschitz_cond} and combining with the above inequality, we have
\begin{equation*}\label{lip_proof_eq1}
   \langle x^*-y^*,y-x\rangle\leq L\norm{x-y}^2.
\end{equation*}
We know that $\norm{x^*-y^*}=\sup_{\norm{u}\leq 1}\langle x^*-y^*,u\rangle$ (see \cite[Section A.1.6]{boyd}). By taking $u=\frac{y-x}{\norm{y-x}}$, it appears from \eqref{lip_proof_eq1} that
\begin{equation}
    \norm{x^*-y^*}\leq \langle x^*-y^*,\frac{y-x}{\norm{y-x}}\rangle \leq L\norm{x-y}.
\end{equation}

To prove part (ii) of this result, let us take $x^*\in \partial_2f(x,x)$ and $y^*\in \partial_2f (y,y)$ arbitrary. Then, by \eqref{subdifferential} and the fact that $f(x,x)=0$ we have
\begin{equation*}
    f(x,y)\geq \langle x^*,y-x\rangle,\quad\forall\, y\in \mathcal{A}.
\end{equation*}
Similarly,
\begin{equation*}
    f(y,x)\geq \langle y^*,x-y\rangle,\quad\forall\, y\in \mathcal{A}.
\end{equation*}
Adding the above two inequalities, we have
\begin{equation*}
    \langle x^*-y^*,y-x\rangle\leq f(x,y)+f(y,x)\leq -\mu\norm{x-y}^2.
\end{equation*}
\end{proof}
\begin{definition}\cite{quoc}
The dual equilibrium problem, denoted by \(\text{EP}_d(f, \mathcal{D})\), aims to find \(\bar{y} \in \mathcal{D}\) such that
\begin{equation}\label{dual_EP}
   f(y, \bar{y}) \leq 0, \quad\forall\, y \in \mathcal{D}.
\end{equation}   
\end{definition}
The solution set of dual problem \eqref{dual_EP} is denoted by \(\text{S}^{d}_{\text{EP}}(f, \mathcal{D})\).

We make the following assumptions on the bifunction $f: X \times X\rightarrow \mathbb{R}$ later in the paper, where the set $X\subset \mathbb{R}^n$ will be specified whenever these assumptions are made. 
\begin{assumption}\label{A1}
The function $f(\cdot,y)$ is upper semi-continuous for all $y \in X$, the function $f (x, \cdot)$ is continuous and convex for all $x \in X$ and $f(x,x)=0$ for all $x\in X$.
\end{assumption}

\begin{lemma}\label{lemma_strng_mono_unique}
\cite[Theorem 1]{blum} Consider a bifunction $f:\mathcal{D}\times \mathcal{D}\to \mathbb{R}$ and $\mathcal{D}\subseteq \mathbb{R}^n$ is non-empty, closed and convex. Under Assumption \textbf{A}\ref{A1}, the solution set of problems EP$(f,\mathcal{D})$ \eqref{EP} and EP$_d(f,\mathcal{D})$ \eqref{dual_EP} coincides if $f$ is monotone. Further, if the bifunction $f$ is strongly monotone, then both problems admit the same unique solution.
\end{lemma}

In the next section, we present the iterative methods to find the projected solution of QEP$(f, \mathcal{K})$ along with their convergence analysis. 
\section{Iterative Algorithms for Projected Solutions}\label{iter_algos}
In this section, we propose two iterative algorithms to find the projected solution of QEP. The first algorithm is developed to address the QEP in which the bifunction $f$ is assumed to be strongly monotone. In the second one, we relax the strongly monotonicity assumption, but as a compensation, we consider nonemptiness of a special type of solution set $S^{*}$ (see \eqref{nonemptyset}). 
\subsection{Algorithm for Strongly Monotone QEP}
To find the projected solution of QEP, we propose an iterative algorithm in which certain convex optimization problems are solved in each iteration. We prove the convergence of this algorithm by assuming that the bifunction $f$ is strongly monotone and the constraint map $\mathcal{K}:\mathcal{C}\rightrightarrows\mathbb{R}^n$ satisfies the following assumption: 
\begin{align*}
    \exists\,\alpha>0~\text{such that}~ \norm{P_{\mathcal{K}(u)} (z)-P_{\mathcal{K}(v)}(z)}\leq \alpha\norm{u-v}~\text{for all}~u,v\in \mathcal{C}~\text{and}~z\in \mathbb{R}^n.
\end{align*}
Note that the above assumption on the map $\mathcal{K}$ is satisfied if $\mathcal{K}$ is a ``moving set" type constraint map, that is, $\mathcal{K}:\mathcal{C}\rightrightarrows\mathbb{R}^n$ is defined as
\begin{equation}\label{map_K}
    \mathcal{K}(x)=\nu(x)+\mathcal{M}_\circ
\end{equation} for a $\alpha$-Lipschitz function $\nu:\mathcal{C}\rightarrow\mathbb{R}^n$ and a non-empty closed convex set $\mathcal{M}_\circ \subset \mathbb{R}^n$ (see \cite[Lemma 3.2]{scrimalli}). 

%\begin{definition}\label{gap_fun}\cite{Gap_fun}
  Let us now recall the concept of dual gap function associated with the dual problem $EP_d(f,\mathcal{D})$ \cite{Gap_fun}.  A function $g:\mathcal{D}\to \mathbb{R}$ is called a gap function for the dual problem \eqref{dual_EP} if $g(y) \ge 0$ for all $y \in \mathcal{D}$, and $g(\bar{y}) = 0$ if and only if $\bar{y} $ solves the dual problem \eqref{dual_EP}.
%\end{definition} 
As per \cite{Gap_fun}, the gap function $g:\mathcal{D}\rightarrow \mathbb{R}^n$ for the dual problem \eqref{dual_EP} satisfying the mentioned conditions is given as follows:
\begin{equation}\label{gap_function}
    g(y) = \sup_{y \in \mathcal{D}}\left\{f(x,y) + \frac{\mu}{2}\|x-y\|^2\right\},
\end{equation}
where $\mu>0$ is the strongly monotone parameter of $f.$ One can observe that $\bar y\in \mathcal{D}$ solves $EP(f,\mathcal{D})$ if $\bar y\in \arg \min_{y\in \mathcal{D}} g(y)$ (see \cite[Remark 2]{quoc} or \cite{Gap_fun}). This aspect of dual gap functions attracted many researchers in the recent past to formulate an algorithm for equilibrium problems in terms of strongly convex minimization problems and find the solution of EP numerically under different sets of assumptions (see e.g.,\cite{scrimalli,quoc,bigi}). In \cite{scrimalli} and \cite{quoc}, the authors have proposed algorithms to find solutions of variational inequalities and equilibrium problems, respectively, using dual gap functions. 

 We first provide an algorithm to solve the equilibrium problem EP($f,\mathcal{D}$) (\ref{EP}) before giving a methodical statement of the first iterative algorithm to find the projected solution of QEP$(f, \mathcal{K})$. The following algorithm to find solutions of EP($f,\mathcal{D}$) will be used in each iteration of our first algorithm to find projected solutions. This algorithm has been adopted from \cite{scrimalli} and \cite{quoc}. %where it is studied to solve variational inequality and equilibrium problems, respectively.

\begin{algorithm}[H]
\caption{Method to solve strongly monotone EP} \label{Algo4} 
\begin{algorithmic}[4] 
\Statex \textbf{Initialization}: Provide $\lambda_0=1,$  a natural number $N$ and  $\hat y\in \mathcal{D}$. Give the values of parameters $L>0$ and $\mu>0$, the Lipschitz constant and coefficient of strong monotonicity, respectively. Solve the strongly convex programming problem
\begin{equation}\label{w_0}
    \arg\min_{w\in \mathcal{D}} \left\{f(\hat y, w)+ \frac{L}{2} \norm{w-\hat y_k}^2\right\}
\end{equation}
to obtain the unique solution $w_{0}$.

\Statex \textbf{Iteration}: For $n=0,\cdots,  N-1$ perform the following three steps:
\Statex \textbf{S1:}  Solve the strongly convex programming problem 
\begin{equation*}
    \min_{z\in \mathcal{D}} \left\{  \sum_{i=0}^n \lambda_i\left[ f(w_i,z)+\frac{\mu}{2} \norm{z-w_i}^2\right]\right\}
\end{equation*}
to obtain a unique solution $z_n$.
\Statex \textbf{S2:} Solve the second strongly convex programming problem
\begin{equation*}
    \min_{w\in \mathcal{D}} \left\{f(z_n,w)+ \frac{L}{2} \norm{w-z_n}^2\right\}
\end{equation*}
to obtain the unique solution $w_{n+1}$.
\Statex \textbf{S3:} Update $\lambda_{n+1}=\frac{\mu}{L} S_n,$ where $S_n=\sum_{i=1}^{n} \lambda_i$.
\Statex \textbf{Output:} Compute the final output $\tilde{y}_{ N} (\mathcal{D},\hat y)= \frac{1}{S_{N}} \sum_{n=0}^{ N} \lambda_i w_i$.
\end{algorithmic}
\end{algorithm}

The efficiency of the solution of EP($f,\mathcal{D}$) obtained by Algorithm \ref{Algo4} depends on the number of iterations ($N$), as shown in the upcoming result. The following result, which is an extension of \cite[Theorem 4.3]{scrimalli} to the case of equilibrium problems, provides the efficiency estimates for Algorithm \ref{Algo4} in terms of distances to the solution of EP($f,\mathcal{D}$).
\begin{proposition}\label{Prop_inequality}
    Suppose $\mathcal{D}$ is non-empty closed convex subset of $\mathbb{R}^n$. Assume that the bifunction $f:\mathcal{D}\times \mathcal{D}\rightarrow \mathbb{R}$ satisfies the condition \textbf{A} \ref{A1}. Additionally, suppose that
    \begin{enumerate}[label=(\alph*)]
        \item\label{th_2a} the bifunction $f$ satisfies Lipschtz type continuity condition \eqref{Lipschitz_cond} on $\mathcal{D}$ with parameter $L>0$;
        \item\label{th_2b} the bifunction $f$ is strongly monotone on $\mathcal{D}$ with parameter $\mu>0$.
        %\item The function $f (\cdot, y):\mathcal{A}\rightarrow \mathbb{R}$ is upper semi-continuous for all $y \in \mathcal{A}$ and $f (x, \cdot):\mathcal{A}\rightarrow \mathbb{R}$ is proper lower semi-continuous and convex for all $x \in \mathcal{A}$.
    \end{enumerate}
    Then, the Algorithm \ref{Algo4} converges to the unique solution of $EP(f,\mathcal{D})$ in the following way:
    \begin{equation}\label{estimate_itera}
        \norm{\tilde{y}_{\widehat N} (\mathcal{D},\hat y)- \bar y(\mathcal{D})}\leq \frac{5L}{\mu} \exp\left\{-\frac{N\mu}{2(L+\mu)}\right\} \norm{\hat y-\bar y(\mathcal{D})},
    \end{equation}
    where $\bar y(D)$ is a unique solution of $EP(f,\mathcal{D})$.
\end{proposition}
\begin{proof}
     We claim that the dual gap function defined in (\ref{gap_function}) satisfies,
    \begin{equation}\label{ineq2}
        g(w_0)\leq \frac{8L^2}{\mu} \norm{\hat y-\bar y(\mathcal{D})}^2.
    \end{equation}
    Since $w_0$ solves the strongly convex programming problem (\ref{w_0}), we have
    \begin{equation}\label{eq_1_strong_EP}
        f(\hat y,w_0)+\frac{L}{2} \norm{\hat y-w_0}^2+\frac{L}{2}\norm{y-w_0}^2\leq f(\hat y,y)+\frac{L}{2} \norm{\hat y-y}^2,\quad\forall\, y\in \mathcal{D}.
    \end{equation}
    By substituting $y=\bar y(D)$ in above inequality, we obtain
    \begin{align}
        \frac{L}{2}\norm{w_0-\hat y}^2\notag &\overset{\ref{th_2b}}\leq -f(\bar y(\mathcal{D}),\hat y)-\mu \norm{\hat y-\bar y(\mathcal{D})}^2-f(\hat y, w_0)+\frac{L}{2}\left\{\norm{\hat y-\bar y(\mathcal{D})}^2-\norm{\bar y(\mathcal{D})-w_0}^2\right\}\notag\\&\overset{\ref{th_2a}}\leq L\norm{\bar y(\mathcal{D})-\hat y}\norm{\hat y-w_0}-\mu \norm{\hat y-\bar y(\mathcal{D})}^2+\frac{L}{2}\left\{\norm{\hat y-\bar y(\mathcal{D})}^2-\norm{\bar y(\mathcal{D})-w_0}^2\right\}.\notag
        %&\leq -\frac{L}{2}\norm{\hat y-w_0}^2+2L\norm{\bar y(\mathcal{D})-\hat y}\norm{\hat y-w_0}-\mu \norm{\hat y-\bar y(\mathcal{D})}^2.\notag
    \end{align}
    \noindent From triangle inequality, we know
    \begin{equation}\label{traingle_ineq}
     \norm{\bar y(\mathcal{D})-w_0}^2\leq \norm{\bar y(\mathcal{D})-\hat y}^2+\norm{\hat y-w_0}^2+ 2\norm{\bar y(\mathcal{D})-\hat y}\norm{\hat y-w_0} .
    \end{equation}
    By using \eqref{traingle_ineq} and the inequality preceding it, we have 
    \begin{equation*}
        {L}\norm{\hat y-w_0}^2+\mu \norm{\hat y-\bar y(\mathcal{D})}^2\leq 2L\norm{\bar y(\mathcal{D})-\hat y}\norm{\hat y-w_0}.
    \end{equation*}
    Hence, it appears that
    \begin{equation}\label{ineq1}
    \frac{\mu}{2L}\norm{\hat y-\bar y(\mathcal{D})}\leq \norm{\hat y-w_0}\leq 2 \norm{\hat y-\bar y(\mathcal{D})}.
    \end{equation}
Now, by definition of dual gap function $g$, we have
\begin{small}
    \begin{eqnarray*}
        g(w_0)%=&\sup_{z\in \mathcal{D}} \left\{f(z,w_0)+\frac{\mu}{2} \norm{z-w_0}^2\right\}\notag\\
       & \overset{\ref{th_2b}}\leq &\sup_{z\in \mathcal{D}} \left\{-f(w_0,z)-\frac{\mu}{2} \norm{z-w_0}^2\right\}\notag\\
      %  &\overset{}\leq &\sup_{z\in \mathcal{D}} \left\{-f(w_0,z)+f(\hat y,z)-f(\hat y,w_0)-\frac{L}{2} \norm{\hat y-w_0}^2+\frac{L}{2} \norm{\hat y-z}^2-\frac{L+\mu}{2} \norm{z-w_0}^2\right\}\notag\\
        &\overset{\eqref{eq_1_strong_EP},\ref{th_2a}}\leq & \sup_{z\in \mathcal{D}} \left\{ L\norm{w_0-z}\norm{\hat y-w_0}-\frac{L}{2} \norm{\hat y-w_0}^2+\frac{L}{2} \norm{\hat y-z}^2-\frac{L+\mu}{2} \norm{z-w_0}^2\right\}\notag\\
        %\leq &\left\{ -\frac{L}{2} \norm{z-w_0}^2-\frac{\mu}{2}\norm{z-w_0}^2\right\}\\
        &\leq& \sup_{z\in \mathcal{D}} \left\{ L\norm{w_0-z}\norm{\hat y-w_0}+L\langle \hat y-w_0,w_0-z\rangle-\frac{\mu}{2} \norm{z-w_0}^2\right \}\notag\\
        &\leq &\frac{2L^2}{\mu} \norm{\hat y-w_0}^2\overset{\eqref{ineq1}}\leq  \frac{8L^2}{\mu} \norm{\hat y-\bar y(\mathcal{D})}^2.\notag
     \end{eqnarray*}
    \end{small}
   From \cite[Theorem 2]{quoc}, we know that the Algorithm \ref{Algo4} converges in the following way:
   %\begin{small}
    \begin{eqnarray*}
        \frac{\mu}{2}\norm{\tilde{y}_{\widehat N} (\mathcal{D},\hat y)- \bar y(\mathcal{D})}&\leq& \left[ g(w_0)+\frac{(L^2-\mu^2)}{2\mu}\norm{w_0-\bar y(\mathcal{D})}^2\right]\exp\left\{-\frac{N\mu}{2(L+\mu)}\right\},\\
        &\overset{\eqref{traingle_ineq},\eqref{ineq1}}{\leq}& \left[ g(w_0)+\frac{9(L^2-\mu^2)}{2\mu}\norm{\hat y-\bar y(\mathcal{D})}^2\right]\exp\left\{-\frac{N\mu}{2(L+\mu)}\right\}.
    \end{eqnarray*}
    %\end{small}
    %By applying \eqref{traingle_ineq} and \eqref{ineq1} in the above inequality, we have
    %\begin{equation*}
     %   \frac{\mu}{2}\norm{\tilde{y}_{\widehat N} (\mathcal{D},\hat y)- \bar y(\mathcal{D})}\leq \left[ g(w_0)+\frac{9(L^2-\mu^2)}{2\mu}\norm{\hat y-\bar y(\mathcal{D})}^2\right]\exp\left\{-\frac{N\mu}{2(L+\mu)}\right\}
    %\end{equation*}
    We observe that the result follows by combining the above inequality with the inequality (\ref{ineq2}). 
\end{proof}
We now propose an algorithm to find the projected solution for QEP($f,\mathcal{K}$). The complexity of the proposed algorithm is relevant to the contraction gap $\delta=1-\alpha \frac{L}{\mu}$. From the estimate (\ref{estimate_itera}), we can now define the minimum number of iterations $\hat N = N(\alpha,\mu,L)$ required to satisfy the condition,
\begin{equation}\label{Remark_estimate}
    \frac{5L}{\mu} \exp\left\{-\frac{\hat N\mu}{2(L+\mu)}\right\} \leq \frac{\delta}{4} \implies \hat N= \left\lceil\frac{2(L+\mu)}{\mu}\ln\left( \frac{20L}{\mu-\alpha L}\right)\right \rceil.
\end{equation}
This estimate allows us to propose the following algorithm for finding the projected solution of QEP($f,\mathcal{K}$).

\begin{algorithm}[H]
%\setstretch{1.2}
\caption{Method to find the projected solution of strongly monotone QEP}
\label{Algo3}
\begin{algorithmic}[1]
\Statex {$\textbf{Initialization:}$ Give an initial point  $y_{0}\in \mathbb{R}^n$ and choose tolerance $\varepsilon>0.$
\Statex \textbf{Step 1}: Calculate $k_{\epsilon}$ from equation \eqref{no_iter}. Set $k  \leftarrow  0.$ }
\Statex\textbf{Iteration:} For $k=0,1,\ldots, k_{\epsilon} $
\Statex \textbf{Step 2}:  Find 
\begin{equation}\label{algo3_eq1}
    x_{k}=P_\mathcal{C}(y_{k})~\text{and}~\hat y_k=P_{\mathcal{K}(x_k)}(y_k).
\end{equation}
\Statex \textbf{Step 3}:  For given \(x_k~\text{and}~  \hat y_k\), apply Algorithm \ref{Algo4} with inputs as $\mathcal{D} = \mathcal{K}(x_{k})$ and $\hat y_k \in \mathcal{K}(x_{k})$. Update $y_{k+1}$ as follows:
\begin{equation*}\label{EP_x}
y_{k+1}= \tilde{y}_{N} (\mathcal{K}(x_k),\hat y_k).
\end{equation*}
\Statex\textbf{Output}: $\bar{x} = P_{\mathcal{C}}(y_{k_{\epsilon}})$, the projected solution of QEP$(f,\mathcal{K})$ with accuracy $\varepsilon.$
\end{algorithmic}
\end{algorithm}
  
% \begin{algorithm}[H]
% \caption{Method to find projected solution of strongly monotone QEP} \label{Algo3}
% \begin{algorithmic}[3] 
% \Statex \textbf{Step 0}: Choose tolerance $\varepsilon>0$ and $y_0\in \mathbb{R}^n$. Set $k \leftarrow  0.$ %Take a bounded auxiliary sequence of positive parameters \(\{ \gamma_k \}\), tolerance $\varepsilon>0$ and choose $y_{0}\in \mathbb{R}^n .$
% \Statex \textbf{Step 1}: (Outer Loop) Find $x_{k}=P_\mathcal{C}(y_{k})$ and $\hat y_k=P_{\mathcal{K}(x_k)}(y_k)$. %Let $w_0=\arg \min_{w\in \mathcal{K}(x_0)} \psi^L_{$
% \Statex \textbf{Step 2:} (Inner Loop)  For given \(x_k~\text{and}~  \hat y_k\), find 
% \begin{equation}\label{EP_x}
% y_{k+1}= \tilde{y}_{\widehat N} (\mathcal{K}(x_k),\hat y_k),
% \end{equation}
% by applying Algorithm \ref{Algo4}.

% %\Statex \textbf{Step 2}: Calculate $x_{k+1}=P_{\mathcal{C}}(y_{k+1}).$
% \Statex \textbf{Step 3}: If $ \|y_k- y_{k+1}\|<\varepsilon,$ Stop. Otherwise, set $k\leftarrow k+1,$ and go to Step 1.
% \end{algorithmic}
% \end{algorithm}
%To show the convergence of Algorithm \ref{algo}, we would first like to derive sufficient conditions for the uniqueness of the projected solutions motivated by \cite[Theorem 4.1]{scrimalli} and \cite[Corollary 2]{scrimalli}, where sufficient conditions for the unique solution of classical QVI are derived. 
For the obtained values of $x_{k}$ and $\hat y_{k}$ in Step 1, we find $\tilde y_{N} (\mathcal{K}(x_k),\hat y_k)$ in Step 2 by taking $\hat y=\hat y_k$ and $\mathcal{D}=\mathcal{K}(x_k)$ in Algorithm \ref{Algo4}. %For this purpose, we adopt the iterative algorithm studied in \cite{scrimalli} and \cite{quoc} to solve the variational inequalities and equilibrium problems, respectively.

%The number of iterations $\hat N$ required to compute $\tilde{y}_{\widehat N} (\mathcal{K}(x_k),\hat y_k)$ which is an approximation for $S_{EP} (f,\mathcal{K}(x_k))$ in Algorithm \ref{Algo4} can be estimated by taking $\mathcal{D}=\mathcal{K}(x_k)$ in the following result (see Remark \ref{number_iterations}). This result extends \cite[Theorem 4.3]{scrimalli} to the case of equilibrium problems.

%\begin{remark}\label{number_iterations} 

%\end{equation}
%\end{remark}
%taken as $\hat N=[2(\gamma+1)\ln{\frac{12 \gamma}{1-\alpha \gamma}-1}]$ motivated by \cite{scrimalli}. In fact, the following result extends 
Suppose the constraint map $\mathcal{K}:\mathcal{C}\rightrightarrows \mathbb{R}^n$ admits non-empty closed convex values and the set $\mathcal{A}=\overline{co}(\mathcal{K}(\mathcal{C})).$ Let $f:\mathcal{A}\times \mathcal{A}\to \mathbb{R}$ satisfies assumption \textbf{A} \ref{A1} and is strongly monotone with parameter $\mu>0$.  We define an operator $T:\mathcal{C}\rightarrow\mathbb{R}^n$ satisfying following relations motivated by the relaxation operator in \cite[Section 4]{scrimalli},
    \begin{equation}\label{T_EP}
    	T(x)\in \mathcal{K}(x)~\text{such that}~f(T(x),z)\geq 0,\quad\forall\, z\in \mathcal{K}(x).
    \end{equation}
     Note that the above operator $T$ is well defined in the view of Lemma \ref{lemma_strng_mono_unique}.
     Let $F:\mathcal{A}\rightarrow\mathbb{R}^n$ is defined as $F(x)=\partial_2 f(x,x)$. Then, it is well known that the solution sets of EP($f,\mathcal{D}$) and VI($F,\mathcal{D}$) coincides if $f$ satisfies assumption \textbf{A} \ref{A1} and $\mathcal{D}$ is closed convex subset of $\mathbb{R}^n$ (see \cite{Aussel_Dutta_Pandit}). Hence, we can say that $T(x)$ satisfies \eqref{T_EP} iff it satisfies the condition \eqref{T3.1_eq1} for any $x\in \mathcal{C}$. Hence, we can say that the operator $T:\mathcal{C}\rightarrow\mathbb{R}^n$ is defined in the following way
    \begin{equation}\label{T3.1_eq1}
        T(x)\in \mathcal{K}(x)~\text{such that}~\exists\, y^*(x)\in F(T(x)),
        \langle y^*(x),z-T(x)\rangle\geq 0,\quad\forall\,z\in \mathcal{K}(x).
    \end{equation} 
    The following result derives the sufficient conditions to ensure the convergence of Algorithm \ref{Algo3} using this relaxation type operator $T$ defined in \eqref{T3.1_eq1}.
\begin{theorem}\label{Thm_strong_mono}
    Suppose $\mathcal{C}$ is closed convex subset of $\mathbb{R}^n$ and $\mathcal{K}:\mathcal{C}\rightrightarrows\mathbb{R}^{n}$ admits closed convex values. Assume that $\mathcal{A}=\overline{co} (\mathcal{K}(\mathcal{C}))$ and $f:\mathcal{A}\times \mathcal{A}\rightarrow \mathbb{R}$ satisfies the assumption \textbf{A} \ref{A1}. %and the map $\mathcal{K}$ be defined as (\ref{map_K}). 
    Additionally, suppose that
    \begin{enumerate}[label=(\alph*)]
        \item \label{hy_3.1i} the bifunction $f$ satisfies Lipschtz type continuity condition \eqref{Lipschitz_cond} on $\mathcal{D}$ with parameter $L>0$;
         \item \label{hy_3.1iii} the bifunction $f$ is a strongly monotone bifunction with parameter $\mu>0$;
        \item \label{hy_3.1ii} there exists $\alpha\in (0,\frac{\mu}{L})$ such that $\|P_{\mathcal{K}(u)} (z)-$$P_{\mathcal{K}(v)}(z)\|$$\leq \alpha\|u-v\|$ for all $u,v\in \mathcal{C}$ and $z\in \mathbb{R}^n$.       
    \end{enumerate}
    Then, there exists a unique projected solution $\bar x$ for QEP$(f,\mathcal{K})$. Further, Algorithm \ref{Algo3} converges to this unique projected solution $\bar x$ in the following way:
    \begin{equation}\label{result_strong_mono_eq}
        \norm{x_k-\bar x}\leq \frac{1}{\delta} \exp \left \{\frac{-\delta}{2}k\right\} \norm{y_0-T(x_0)},
    \end{equation}
    where $\delta=1-\alpha\frac{L}{\mu}$.
\end{theorem}
\begin{proof}
Let us define the operator $T:\mathcal{C}\rightarrow\mathbb{R}^n$ as \eqref{T3.1_eq1}. It can be observed that the solution set of EP($f,\mathcal{K}(x)$) coincides with the solution set of VI($F,\mathcal{K}(x)$) where $F(x)=\partial_2 f(x,x)$ (see \cite{Aussel_Dutta_Pandit}). Hence, it follows that $x$ is a projected solution of $QEP(f,\mathcal{K})$ if and only if $x$ is a fixed point of $P_\mathcal{C}\circ T$. We claim that the map $P_\mathcal{C}\circ T:\mathcal{C}\rightarrow \mathcal{C}$ is a contraction. 

For any $x\in \mathcal{C}$ and $\lambda>0$, we have 
\begin{equation*}
T(x)\overset{\eqref{projection_inequality},\eqref{T3.1_eq1}}{=} P_{\mathcal{K}(x)} (T(x)-\lambda y^*(x))~\text{for some}~y^*(x)\in \partial_2 f(T(x),T(x)).
\end{equation*}
Fix $x_1,x_2\in \mathcal{C}$. Suppose $z_2=P_{\mathcal{K}(x_2)} (T(x_1)-\lambda y^*(x_1))$. Then, we have
\begin{equation}\label{thm_st_mono_eq1}
    \norm{T(x_1)-z_2} \overset{\ref{hy_3.1ii}}\leq \alpha \norm{x_1-x_2}.
\end{equation}
Since $z_2=Pr_{\mathcal{K}(x_2)} (T(x_1)-\lambda y^*(x_1))$, we observe
\begin{small}
\begin{align*}
    \langle z_2-T(x_1),T(x_2)-z_2\rangle&~~\geq \lambda\langle y^*(x_1),z_2-T(x_2)\rangle\\
    &~~=\lambda \langle y^*(x_1),z_2-T(x_1)\rangle+\lambda\langle y^*(x_1),T(x_1)-T(x_2)\rangle\\
%&\overset{}{\geq} \lambda [\langle y^*(x_1),z_2-y(x_1)\rangle+\langle y^*(x_2),y(x_1)-y(x_2)\rangle+\mu \norm{y(x_1)-y(x_2)}^2]\\
    %&~~\geq \lambda\langle y^*(x_1)-y^*(x_2),z_2-y(x_1)\rangle +\lambda\langle y^*(x_2),z_2-y(x_2)\rangle+ \lambda \mu \norm{y(x_1)-y(x_2)}^2\\
    &\overset{\textbf{L}\ref{Lemma_strong_mono},(\ref{T3.1_eq1})}{\geq}\lambda\langle y^*(x_1)-y^*(x_2),z_2-T(x_1)\rangle +\lambda \mu \norm{T(x_1)-T(x_2)}^2.
\end{align*}
\end{small}
Since $\langle z_2-T(x_1),T(x_2)-T(x_1)\rangle=\langle z_2-T(x_1),T(x_2)-z_2\rangle+\norm{z_2-T(x_1)}^2$, it appears from above inequality that
\begin{align*}
   \lambda \mu \norm{T(x_1)-T(x_2)}^2&~\,\leq \lambda \langle y^*(x_1)-y^*(x_2),T(x_1)-z_2\rangle+\langle z_2-T(x_1),T(x_2)-T(x_1)\rangle\\
   &\overset{\ref{hy_3.1i},\textbf{L}\ref{Lemma_strong_mono}}\leq (1+\lambda L) \norm{T(x_1)-T(x_2)}\norm{z_2-T(x_1)}. 
\end{align*}
Combining the above inequality with the inequality (\ref{thm_st_mono_eq1}), we have 
\begin{equation*}
  \lambda \mu \norm{T(x_1)-T(x_2)} \leq \alpha (1+\lambda L) \norm{x_1-x_2}.  
\end{equation*}
Since $\lambda>0$ can be chosen sufficiently large, we have 
\begin{equation}\label{Thm_str_mono_eq2}
    \norm{T(x_1)-T(x_2)} \leq \frac{\alpha L}{\mu}\norm{x_1-x_2}.  
\end{equation}
Since $P_{\mathcal{C}}$ is non-expansive, we have
\begin{align*}
   \norm{ P_{\mathcal{C}}(T(x_1))-P_{\mathcal{C}}(T(x_2))}&\leq \norm{T(x_1)-T(x_2)}\overset{\eqref{Thm_str_mono_eq2}}{\leq} \frac{\alpha L}{\mu}\norm{x_1-x_2}.
\end{align*}
Hence $P_{\mathcal{C}}\circ T$ satisfies the Banach contraction principle according to the hypothesis \ref{hy_3.1ii} and admits a unique fixed point $\bar x\in \mathcal{C}$, which is a unique projected solution for QEP($f,\mathcal{K}$). 

In the next part of the proof, we show that the sequence $(x_k)_{k\in \mathbb{N}}$ generated by Algorithm \ref{Algo3} converges to $\bar x$. Let us consider $\norm{y_k-T(x_k)}=r_k$.
%\begin{equation}\label{T_3.1_eq6}
%\end{equation}
%Since $\hat y_k=P_{\mathcal{K}(x_k)}(y_k)$. 
In view of Lemma \ref{Lemma_strong_mono}, the assumptions in Proposition \ref{Prop_inequality} are satisfied. Hence, by using inequalities \eqref{estimate_itera} and (\ref{Remark_estimate}), we obtain
 \begin{eqnarray}
     \norm{y_{k+1}-T(x_k)}&\leq& \frac{\delta}{4} \norm{\hat y_k-T(x_k)}\notag\\ &\overset{(\ref{algo3_eq1})}{\leq}& \frac{\delta}{4} \norm{y_k-T(x_k)}.
\label{Thm_str_mono_eq3} \end{eqnarray}
 Also, we have 
 \begin{equation}\label{Thm_str_mon_eq4}
     \norm{T(x_{k+1})-T(x_k)}\overset{\eqref{Thm_str_mono_eq2}}{\leq} \frac{\alpha L}{\mu} \norm{x_{k+1}-x_k}\overset{(\ref{algo3_eq1})}{\leq}\frac{\alpha L}{\mu} \norm{y_{k+1}-y_k}.
 \end{equation}
 Therefore, we have 
 \begin{eqnarray}
   r_{k+1}&=&\norm{y_{k+1}-T(x_{k+1})}{\leq} \norm{y_{k+1}-T(x_k)}+\norm{T(x_{k+1})-T(x_k)}\notag\\&\overset{\eqref{Thm_str_mono_eq3},\eqref{Thm_str_mon_eq4}}{\leq}&
     \frac{\delta}{4} \norm{y_k-T(x_k)}+\frac{\alpha L}{\mu} \norm{y_{k+1}-y_k}  \notag\\&\leq &\left( \frac{\alpha L}{\mu}+\frac{\delta}{4}  \right) r_k+ \frac{\alpha L}{\mu} \norm{y_{k+1}-T(x_k)} \notag
\\&\overset{\eqref{Thm_str_mono_eq3}}\leq &\left(1-\frac{\delta}{2}\right) r_k\leq \exp\left\{ -\frac{\delta}{2}(k+1)\right\}r_0. \notag  
 \end{eqnarray}
 Thus, we have
 \begin{equation}\label{Thm_str_mono_eq7}
     r_k\leq \exp\left\{ \frac{-\delta }{2}k\right\}r_0.
 \end{equation}
 From (\ref{Thm_str_mono_eq2}) and the fact that $x_k=P_{\mathcal{C}} (y_k)$, we notice
 \begin{align}\label{Thm_str_mono_eq8}
     r_k=\norm{y_k-T(x_k)}\overset{\eqref{Thm_str_mono_eq2}}{\geq} \norm{y_k-T(\bar x)}- \frac{\alpha L}{\mu} \norm{ \bar x- x_k}\geq \delta \norm{\bar x-x_k}. 
 \end{align}
We obtain the inequality \eqref{result_strong_mono_eq} by combining the inequalities (\ref{Thm_str_mono_eq7}) and (\ref{Thm_str_mono_eq8}). Hence, the convergence of $(x_k)_{k\in \mathbb{N}}$ to $\bar x$ follows.
 %For this purpose, we first claim that 
 %\begin{equation}
  %   \norm{\tilde y_{\widehat N}(\mathcal{K}(x_k),\hat y_k)-y(x_k)}\leq 3\gamma \exp\left\{ -\frac{\widehat N}{2(\gamma+1)}\right\}\norm{\hat y-y(x_k)} 
 %\end{equation} 
 %for any $k\geq 1$. In fact, we have 
 %\begin{align}
  %   \norm{\tilde y_{\widehat N}(\mathcal{K}(x_k),\hat y_k)-y(x_k)}\leq \frac{2}{\mu} [g(x_\circ)+ ]
 %\end{align}
 %Suppose that the tolerance $\delta_\circ$ for the inner loop is set to $\delta_\circ= \frac{\delta}{4} \norm{\hat y_k-y(x_k)}$ where $\delta=1-\alpha \gamma$. We claim that for any $k\geq 0$ the inner loop converges in atmost $\widehat N=[2(\gamma+1)\ln{\frac{12 \gamma}{1-\alpha \gamma}-1}]$ steps.
\end{proof}
The following result provides the estimation of iterations required in Algorithm \ref{Algo3} to compute the projected solution of QEP($f,\mathcal{K}$) with desired accuracy. This result is a consequence of Proposition \ref{Prop_inequality} and Theorem \ref{Thm_strong_mono}.
\begin{corollary} \label{Coro_str_mono}
    Suppose $\gamma=\frac{L}{\mu}.$ In order to obtain an approximation $x_k$ of exact projected solution $\bar x$ for QEP($f,\mathcal{K}$) with the accuracy $\norm{x_k-\bar x}\leq \epsilon$, the Algorithm \ref{Algo3} requires at most $k_\epsilon$ iterations, given as
    \begin{equation}\label{no_iter}
        k_\epsilon=\left \lceil\frac{2}{1-\alpha \gamma} \ln \frac{\norm{y_0-T(x_0)}}{\epsilon(1-\alpha\gamma)}\right\rceil.
    \end{equation}
    Hence, the total number of computations of bifunction $f(x,y)$ are given as,
    \begin{equation*}
        \left \lceil \frac{10(1+\gamma)}{1-\alpha \gamma}\ln\left( \frac{20\gamma}{1-\alpha \gamma}\right) \ln \frac{\norm{y_0-T(x_0)}}{\epsilon(1-\alpha\gamma)}\right\rceil.
    \end{equation*}
\end{corollary}
\begin{remark}
     Note that \cite[Theorem 4.4]{scrimalli} and \cite[Corollary 3]{scrimalli} show the convergence of an iterative algorithm similar to Algorithm \ref{Algo3} to a classical solution of QVI problems. If the constraint map $\mathcal{K}:\mathbb{R}^n\rightrightarrows \mathbb{R}^n$, Theorem \ref{Thm_strong_mono} and Corollary \ref{Coro_str_mono} assure the convergence of Algorithm \ref{Algo3} to a classical solution of QEP. Hence, Theorem \ref{Thm_strong_mono} and Corollary \ref{Coro_str_mono} extend \cite[Theorem 4.4]{scrimalli} and \cite[Corollary 3]{scrimalli} to the case of QEP.  
\end{remark}
%\begin{theorem}
%     Let the function $f$ be differentiable and let the map $\mathcal{K}$ be defined as (\ref{map_K}). Then Algorithm \ref{algo} converges to a unique projected solution if $k<\frac{\mu}{L}$.
%\end{theorem}

%\begin{remark}
  
%\end{remark}

\subsection{ Proximal Point Method and its Convergence Analysis}
In this section, a proximal point method is presented in which the projected solution for QEP$(f,\mathcal{K})$ is obtained by iteratively solving equilibrium problems.  We show that the sequence generated by this algorithm converges under the following assumptions on the bifunction $f: \mathcal{A} \times \mathcal{A} \to \mathbb{R}$ and the map $\mathcal{K}:\mathcal{C}\rightrightarrows \mathbb{R}^n$:
\begin{assumption}\label{A2}
    The bifunction $f$ is monotone.
\end{assumption}
\begin{assumption}\label{A3}
		%\item[\textbf{(A1)}] \(f(x, x) = 0\) for all \(x \in \mathcal{C}\);
		%\item[\textbf{(A2)}] \(f( \cdot, \cdot ) : \mathcal{A} \times \mathcal{A} \to \mathbb{R}\) is jointly continuous on \(\mathcal{A} \times \mathcal{A}\);
		%\item[\textbf{(A3)}] \(f(x, \cdot) : \mathcal{A} \to \mathbb{R}\) is convex for all \(x \in \mathcal{C}\);
    %\item[\textbf{(A4)}] \textbf{Monotone bifunction:} $f(x, y) + f(y, x) \leq 0$ for all $x, y \in \mathcal{A}$;
     For any sequence $(y_n)_{n \in \mathbb{N}}\subset \mathcal{A}$ with $\lim_n\norm{y_n}=\infty$, there exists $z\in \mathcal{A}$ and $n_0\in \mathbb{N}$ such that $f(y_n,z)\leq 0$ for all $n\geq n_0$.
\end{assumption}

Further, we need some continuity assumptions on the map $\mathcal{K}$ (see for e.g., \cite{QEP_proximal}) for quasi-equilibrium problems. Let us recall the concept of closed and lower semi-continuous maps. A map $\mathcal{K}$ is said to be:
\begin{itemize}
    \item \textit{closed} if for any $(x_k)_{k\in \mathbb{N}} \subset \mathcal{C}$ with $x_k \to x$ and $y_k \in \mathcal{K}(x_k)$ with $y_k \to y$, we have that $y \in \mathcal{K}(x)$;
    \item \textit{lower semi-continuous} if for any sequence $(x_k)_{k\in \mathbb{N}} \subset \mathcal{C}$ with $x_k \to x$ and for each $y \in \mathcal{K}(x)$ there exists a sequence $(y_k)_{k\in \mathbb{N}} \subset \mathcal{A}$ with $y_k \in \mathcal{K}(x_k)$ such that $y_k \to y$.
\end{itemize}

\begin{assumption}\label{A4}
     $\mathcal{K}$ is a closed lower semi-continuous map with non-empty convex closed values. 
\end{assumption}

%When the strict inequality holds in (A4), for $y \neq x$, $f$ is called strictly monotone. 
Let us recall the following result from \cite{blum,sosa_EP} on the non-emptiness of the solution set for equilibrium problems.

\begin{lemma}\cite{blum}
Assume that the bifunction $f:\mathcal{A}\times \mathcal{A}\rightarrow \mathbb{R}$ satisfies the assumptions \textbf{A} \ref{A1}-\textbf{A} \ref{A3}. Then, $\text{S}_{\text{EP}}(f, \mathcal{B}) = \text{S}^{d}_{\text{EP}}(f, \mathcal{B})\neq \emptyset$ for any non-empty closed convex set $\mathcal{B}\subseteq \mathcal{A}$.
\end{lemma}

%\begin{proposition}
%Assume that the bifunction $f:\mathcal{A}\times \mathcal{A}\rightarrow \mathbb{R}$ satisfies the assumptions \textbf{A} \ref{A1}-\textbf{A} \ref{A3}. Then $\text{S}_{\text{EP}}(f, \mathcal{B})$ is nonempty for any nonempty closed convex set $\mathcal{B}\subseteq \mathcal{A}$.
%\end{proposition}
The regularization of bifunctions plays an important role in obtaining the unique solution of EP at each iteration in the proximal point method (see \cite{sosa_EP,QEP_proximal}). Let us recall the concept of regularization of $f$ from \cite{sosa_EP}. %Before presenting the method, let us state a useful result on regularization methods for EP, whose proof follows from \cite[Proposition 4]{sosa_EP}, by letting $K = \mathcal{K}(\bar{x})$ and taking into account that $\bar{y}\in \mathcal{A}$ is fixed.
For some fixed $\gamma > 0$ and $\bar{y}\in \mathbb{R}^n$, a bifunction $\hat f:\mathcal{A}\times \mathcal{A}\to \mathbb{R}$ defined as
\begin{equation}\label{regularized_f}
  \hat{f}(x, y) = f(x, y) + \gamma \langle x - \bar{y}, y - x \rangle.  
\end{equation}
is known as \textit{regularization} of $f$. Further, the problem EP$(\hat f,\mathcal{B})$ for some $\mathcal{B}\subseteq \mathcal{A}$ is known as the regularized equilibrium problem.
The following result obtained by combining \cite[Proposition 1]{sosa_EP} and \cite[Proposition 3]{sosa_EP} shows the advantage of regularizing the bifunction $f$ in (\ref{regularized_f}).
\begin{lemma}\cite{sosa_EP}\label{prop_existence}
    Suppose that the bifunction $f:\mathcal{A}\times \mathcal{A}\rightarrow \mathbb{R}$ satisfies the assumptions \textbf{A} \ref{A1} and \textbf{A} \ref{A2}. Let $\mathcal{B}$ be a non-empty closed convex subset of $\mathcal{A}$. Then $EP(\hat f,\mathcal{B})$ admits a unique solution, where $\hat f$ is given by \eqref{regularized_f}.
\end{lemma}
We recall a result from \cite{sosa_EP}, which shows the relation between the solutions of a dual equilibrium problem and a regularized equilibrium problem.
\begin{lemma}\cite[Proposition 4]{sosa_EP} \label{Prop_Fejer_inequality}
Suppose that the bifunction $f:\mathcal{A}\times \mathcal{A}\rightarrow \mathbb{R}$ satisfies an assumption \textbf{A} \ref{A1}. Let $\mathcal{B}$ be a non-empty closed convex subset of $\mathcal{A}$. %Let $\bar{x} \in \mathcal{C}$ and $\bar{y}\in \mathcal{A}$ be arbitrary points. 
Consider $\hat{y} \in S_{EP}(\hat{f}, \mathcal{B})$ and $y^* \in S^{d}_{EP}(f, \mathcal{B})$, where $\hat{f}$ is given by (\ref{regularized_f}). Then, %If $f:\mathcal{A}\times \mathcal{A}\rightarrow \mathbb{R}$ satisfies \textbf{A} \ref{A1}-\textbf{A} \ref{A3} and $\mathcal{K}$ admits non-empty convex closed values, then
\[
\|\hat{y} - y^*\|^2 + \|\bar{y} - \hat{y}\|^2 \leq \|\bar{y} - y^*\|^2.
\]
\end{lemma}

Let us now recall the concept of Fejér convergence of a sequence to a set. We say that $(y_k)_{k\in \mathbb{N}}\subset \mathbb{R}^n$ is \textit{Fejér convergent} to a nonempty subset $\mathcal{D}$ of $\mathbb{R}^n$ if
\[
\|y_{k+1} - y\| \leq \|y_k - y\|, \quad\forall\, k \in \mathbb{N}, \enspace \forall\,y \in \mathcal{D}.
\]
The notion of Fejér convergence is frequently used to prove the convergence of sequences generated by several algorithms. The following result related to Fejér convergence of sequences is proved by Alber et al. \cite{alber1998projected}.

\begin{lemma}\label{lemma_Fejer_conv}
    Let $\mathcal{D}\subseteq \mathbb{R}^n$ be non-empty and the sequence $(y_k)_{k\in \mathbb{N}}$ be Fejér convergent to $\mathcal{D}$. Then $(y_k)_{k\in \mathbb{N}}$ is a bounded sequence. Furthermore, $(y_k)_{k\in \mathbb{N}}$ converges to $y$ if $y\in \mathcal{D}$ is a cluster point of $(y_k)_{k\in \mathbb{N}}$.
\end{lemma}

% In the view of Lemma \ref{prop_existence}, we give the methodical statement of the algorithm as below:

% Now, we present a proximal point method for solving a QEP$(f, \mathcal{K})$ as well as its convergence analysis, assuming from now on that f satisfies P1–P5 and the set S∗ given by (7) is nonempty.

 Iusem-Sosa \cite{sosa_EP} studied the proximal point method for equilibrium problems. Subsequently, Santos-Souza \cite{QEP_proximal} extended the proximal point method to find the classical solution of QEP$(f, \mathcal{K})$ when $\mathcal{K}(\mathcal{C})\subseteq \mathcal{C}$. Motivated by these approaches, we now describe the proximal point algorithm to find the projected solution for QEP$(f, \mathcal{K})$ when $\mathcal{K}$ need not be a self-map.%The following pseudocode outlines the steps involved in the iterative procedure.

% assumptions \textbf{A} \ref{A1}-\textbf{A} \ref{A3} and $\mathcal{K}:\mathcal{C}\rightrightarrows\mathbb{R}^{n}$ satisfies the assumption \textbf{A} \ref{A4}. Let $S^{*}$ defined in \eqref{nonemptyset} be nonempty.

\begin{algorithm}[H]
\caption{Proximal-point method to find projected solution of QEP} \label{algo} 
\begin{algorithmic}[1] 
\Statex \textbf{Initialization}: Take a bounded auxiliary sequence of positive parameters \(\{ \gamma_k \}\), tolerance $\varepsilon>0$ and choose $y_{0}\in \mathbb{R}^n .$
\Statex \textbf{Step 1}: (Outer loop) Set $x_0=P_{\mathcal{C}}(y_0).$
\Statex \textbf{Step 2:} (Inner loop)  For a given \(x_k,~k \ge 0\), solve the equilibrium problem $EP(f_{k}, \mathcal{K}(x_{k})).$ Let
\begin{equation*}\label{EP_x}
\{y_{k+1}\} = \text{S}_{\text{EP}}(f_k, \mathcal{K}(x_{k})),
\end{equation*}
where \(f_k(x, y) = f(x, y) + \gamma_k \langle x-y_{k},y-x\rangle\).
%\Statex \textbf{Step 2}: Calculate $x_{k+1}=P_{\mathcal{C}}(y_{k+1}).$
\Statex \textbf{Step 3}: If $ \|y_k- y_{k+1}\|<\varepsilon,$ Stop. Otherwise, set $k\leftarrow k+1,$ and go to Step 1 by taking $x_{k+1}=P_\mathcal{C}(y_{k+1})$.
\end{algorithmic}
\end{algorithm}

We observe that the inner loop of Algorithm \ref{algo} (see Step 2) is well defined due to Lemma \ref{prop_existence}. We solve $EP(f_k, \mathcal{K}(x_k))$ in the inner loop for a given value of $x_{k}$ and $y_{k}$ from the outer loop in Step 1. We adopt the inexact subgradient method studied in \cite{santos2011inexact} to solve this equilibrium problem with $\xi_{k}=0$ and $\varepsilon_{k} = 0$ for all $k$.   

\begin{algorithm}[H]
\caption{Inexact projected subgradient method to solve $EP(f_{k},\mathcal{K}(x_{k}))$} \label{algo1} 
\begin{algorithmic}[1] 
\Statex \textbf{S1:} Choose $w_{0}\in  \mathcal{K}(x_{k})$ and tolerance $\delta>0.$ Set $j=0.$
\Statex \textbf{S2:} Obtain $g^j \in \partial_{2}f(w^j,w^j).$ Calculate 
$$\alpha^{j} = \frac{\beta^{j}}{\gamma^{j}},~\text{where } \gamma^{j} = \max\{\rho^{j},\|g^{j}\|\}.$$
\Statex \textbf{S3:} Find $w^{j+1}$ such that $w^{j+1} = P_{\mathcal{K}(x_{k})}(w^j - \alpha^{j}g^{j} )$ 
\Statex \textbf{S4:} If $ \|w^{j+1}-w^j\|<\delta,$ Stop and set $y_{k+1} = w^{j+1}$. Otherwise, set $j\leftarrow j+1,$ and go to S1.
\end{algorithmic}
\end{algorithm}

\begin{figure}[H]
    \centering
    \includegraphics[scale=1.2]{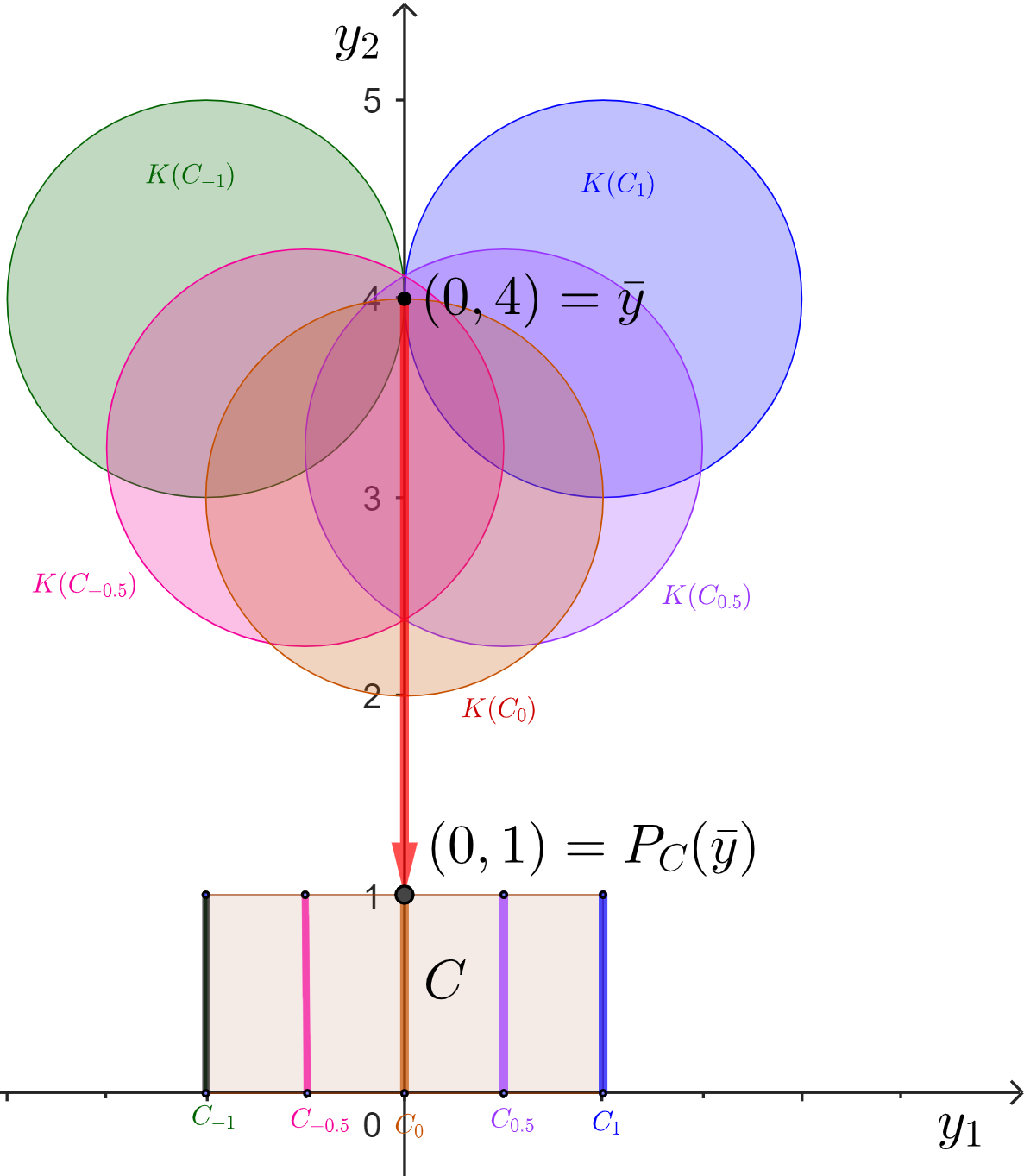}
    \caption{Illustration of projected solution in Example \ref{PPM_ex}}
    \label{PPM}
\end{figure}

 Now, we provide the convergence criteria of Algorithm \ref{algo} based on a special type of solution set, which is defined as follows:
\begin{equation}\label{nonemptyset}
S^* = \left\{ \bar{y} \in \bigcap_{z \in \mathcal{C}} \mathcal{K}(z)\bigg{|}\, f(\bar{y}, y) \geq 0, \quad\forall\, y \in \bigcup_{z \in \mathcal{C}} \mathcal{K}(z) \right\}. 
\end{equation}
Denote 
\begin{equation}\label{S_QEP}
	\tilde S_{QEP}(f,\mathcal{K})=\{\bar y\in \mathcal{A}\,|\,P_{\mathcal{C}}(\bar y)\in S^P_{QEP}(f,\mathcal{K})\}. 
\end{equation} Clearly, $S^* \subset \tilde S_{QEP}(f,\mathcal{K}).$ To prove the convergence of Algorithm \ref{algo}, we assume that $S^* \neq \emptyset.$ Further, we define the set 
\begin{equation*}
    S^{*}_{P} = P_{\mathcal{C}}(S^*) = \bigcup_{\bar{y}\in S^*}P_{\mathcal{C}}(\bar{y}).
\end{equation*}
This assumption was considered to study the convergence of various methods for solving QEP and QVIP, such as proximal point algorithms, extragradient algorithms and projection-like methods. In the following example, we show that there exists a non-empty set $S^*$ which further verifies $S^{*}_{P}$ is nonempty.
\begin{example}\label{PPM_ex}
 Let $QEP(f, \mathcal{K})$ be defined as 
 $\mathcal{C} = [-1,1]\times [0, 1],~~\mathcal{K}(x_1,x_2) = \{(y_1,y_2)\in\mathbb{R}^2: (y_1 - x_1)^2 + (y_2 - 3-x_1^2)^2 \le 1\}$ and $f(x,y) = x_1(y_2 - x_2) + y_1(y_1-x_1).$ Note that for this example, $\mathcal{K}(\mathcal{C})  \cap \mathcal{C} = \emptyset.$ Therefore, there is no classical solution to the QEP$(f, \mathcal{K})$. We can visualize in Figure \ref{PPM} that $\mathcal{K}(\mathcal{C})$ is the union of circles that correspond to a vertical strip in the set $\mathcal{C}.$  Also, there is only one point $(0,4)$ that belongs to $\bigcap\limits_{z \in \mathcal{C}} \mathcal{K}(z).$ Furthermore, $f((0,4),y) = y_1^2 \ge 0,$ for all $y\in \bigcup\limits_{z \in \mathcal{C}} \mathcal{K}(z).$ Hence, $(0,4) \in S^*$. The projection of the point $(0,4)$ is $(0,1)$, which is shown by the red arrow in Figure \ref{PPM}. Therefore, $(0,1) \in S_{P}^{*}$, which is a projected solution of the mentioned $QEP(f, \mathcal{K})$.
\end{example}

In Step 3 of Algorithm \ref{algo}, we used the stopping criteria of the form $\| y_k-y_{k+1}\|<\varepsilon$ to determine when the iterative process should terminate.
This condition helps to check whether the values are getting close enough to each other, which means that the solution is nearly reached. In the next proposition, we show that when $y_{k+1}=y_{k}$, Algorithm \ref{algo} gives the exact solution. This supports that our stopping condition makes sense and works well in Algorithm \ref{algo}.

% In the next proposition, we show that the stopping condition used in Step 3 of Algorithm \ref{algo} of a condition under which the sequence generated by converges to the projected solution.
\begin{proposition}
Let $(y_k)_{k\in \mathbb{N}}$ be the sequence generated by Algorithm \ref{algo} such that $ y_{k+1} = y_k$ for some \(k \geq 0\), then $P_{\mathcal{C}}(y_k)$ is a projected solution of QEP$(f,\mathcal{K})$. 
\end{proposition}

\begin{proof}
From the definition of the Algorithm \ref{algo}, we have that \(y_{k} \in  \text{S}_{\text{EP}}(f_k, \mathcal{K}(x_{k}))\). Thus, \(y_{k} \in \mathcal{K}(x_{k})\), and
$$f_{k}(y_{k+1},y) = f(y_{k+1},y) +\gamma_k \langle y_{k+1} - y_{k}, y - y_{k+1} \rangle \ge 0 \quad\forall\, y \in \mathcal{K}(x_{k}).$$
Since, \(y_{k+1} \in \mathcal{K}(x_{k})\) with $y_{k+1}=y_{k},$ we have $f(y_{k+1},y) \geq 0 $ for all $y \in \mathcal{K}(x_{k}).$ Further, $P_{\mathcal{C}}(y_k)=P_{\mathcal{C}}(y_{k+1})$ implies $x_{k}= P_{\mathcal{C}}(y_{k+1})$. Hence, $x_k$ is a projected solution of QEP.
\end{proof}

In the next part of this section, we show the convergence of the sequence generated by Algorithm \ref{algo}. Before establishing the main convergence result, we first show that the sequence generated by Algorithm \ref{algo} satisfies certain key properties, including Fejér convergence, boundedness, and asymptotic regularity.

\begin{proposition}\label{Prop_Fejer_soln}
 Suppose $f:\mathcal{A}\times \mathcal{A}\rightarrow \mathbb{R}$ satisfies the assumptions \textbf{A} \ref{A1}-\textbf{A} \ref{A3} and $\mathcal{K}:\mathcal{C}\rightrightarrows\mathbb{R}^{n}$ satisfies \textbf{A} \ref{A4}. Let $S^{*}$ defined in \eqref{nonemptyset} be non-empty.  Then the following assertions hold for the sequence $(y_k)_{k\in \mathbb{N}}$ produced by Algorithm \ref{algo}: 
\begin{enumerate}
    \item [(i)]$(y_k)_{k\in \mathbb{N}}$ is Fejér convergent to $S^*$;
    \item [(ii)] $(y_k)_{k\in \mathbb{N}}$ is bounded;
    \item [(iii)] $\lim_{k \to \infty} \|y_{k+1} - y_k\| = 0$.
\end{enumerate}    
\end{proposition}

\begin{proof}
Let us take an arbitrary point $y^*\in S^*\subseteq S_{EP} (f,\mathcal{K}(x))$ for all $x\in \mathcal{D}$. The monotonicity of $f$ implies $y^*\in S^d_{EP} (f,\mathcal{K}(x_k))$ for all $k\in \mathbb{N}$. We have $y_{k+1}\in S_{EP}(f_k,\mathcal{K}(x_k))$ as per Algorithm \ref{algo}. Hence, by taking \( \hat{f} = f_k \), \( \bar{y} = y_k \), \( \hat{y} = y_{k+1} \)in Lemma \ref{Prop_Fejer_inequality}, we obtain
\begin{equation*}
     \| y_{k+1} - y^* \|^2 + \| y_{k+1} - y_k \|^2 \leq \| y_k - y^* \|^2, \quad\forall\, k \in \mathbb{N}.
\end{equation*}

Since $y^*$ was arbitrarily taken in $S^*$, we have $(y_k)_{k\in \mathbb{N}}$ is Fejér convergent to $S^*$. Finally, the assertion (ii) and (iii) follows from Lemma \ref{lemma_Fejer_conv}. 

%Similarly, if \( x^* \in P_\mathcal{C}(S^d_{EP}(\hat{f}, K(\bar{x}))) \), then there exists \( y^* \in S^d_{EP}(\hat{f}, K(\bar{x})) \) such that \( P_C(y^*) = x^* \).

%Hence, the following inequality holds due to Proposition \ref{Prop_dual}:
\end{proof}

Now, we prove the sequence generated by Algorithm \ref{algo} converges to a projected solution of QEP$(f,\mathcal{K})$ by using the key properties of this sequence derived in the above results.

\begin{theorem}\label{Thm_proximal_conv}
    Suppose $f:\mathcal{A}\times \mathcal{A}\rightarrow \mathbb{R}$ satisfies the assumptions \textbf{A} \ref{A1}-\textbf{A} \ref{A3} and $\mathcal{K}:\mathcal{C}\rightrightarrows\mathbb{R}^{n}$ satisfies the assumption \textbf{A} \ref{A4}. Let $S^{*}$ defined in \eqref{nonemptyset} be nonempty. Then, every cluster point $\bar y$ of the sequence $(y_{k})_{k\in \mathbb{N}}$ generated by Algorithm \ref{algo} yields the projected solution $\bar x=P_\mathcal{C}(\bar y)$ of QEP$(f,\mathcal{K})$.
\end{theorem}
\begin{proof}
    Let $(y_{k_{j}})_{{j}\in \mathbb{N}}$ be the subsequence of $(y_{k})_{k\in \mathbb{N}}$ converging to $\bar{y}.$ %Then, there exists $\bar y\in \mathcal{A}$ such that $P_{\mathcal{C}}(\bar y)= \bar x$. Further,
    By Step 1 of Algorithm \ref{algo}, we have $y_{k_{j}+1} \in \mathcal{K}(x_{k_{j}})$, for $x_{k_j}=P_{\mathcal{C}}(y_{k_j})$, such that  $f_{k_{j}}(y_{k_{j}+1},y)\geq 0~\forall~y \in \mathcal{K}(x_{k_{j}}).$ In the view of Proposition \ref{Prop_Fejer_soln}(iii), $y_{k_{j}+1}$ converges to $\bar{y}$. Since $P_{\mathcal{C}}$ is continuous, we observe $P_{\mathcal{C}}(y_{k_{j}})=x_{k_{j}}$ converges to $P_{\mathcal{C}}(\bar{y})=\bar{x}$. Now, we want to show that $\bar{x}$ is the projected solution. It is sufficient to show $f(\bar{y},y)\ge 0$ for all $y \in \mathcal{K}(\bar{x}).$ Let $z \in \mathcal{K}(\bar{x})$ be arbitrary, then there exists $z_{k_{j}}$ converging to $z$ and $z_{k_{j}} \in \mathcal{K}(x_{k_{j}}).$ Since $f_{k_{j}}(y_{k_{j}+1},y)\geq 0$ for all $y\in \mathcal{K}(x_{k_{j}}),$ we have 
    \begin{equation*}
f(y_{k_{j}+1},z_{k_{j}})+\gamma_{k_{j}}\langle y_{k_{j}+1}-y_{k_{j}},z_{k_j}- y_{k_{j}+1}\rangle \geq 0.
    \end{equation*}
    This shows,
     \begin{equation*}
f(y_{k_{j}+1},z_{k_{j}})+\gamma_{k_{j}}\norm{y_{k_{j}+1}-y_{k_{j}}}\norm{z_{k_j}- y_{k_{j}+1}} \geq 0.
 \end{equation*}
   Since, $f$ fulfills assumption \textbf{A} \ref{A1} and $\gamma_k$ is bounded, we have $f(\bar{y},z)\geq 0.$ Therefore, $\bar{y}$ is projected solution.
\end{proof}

We note that Theorem \ref{Thm_proximal_conv} does not guarantee that the cluster point of $(y_k)_{k\in \mathbb{N}}$ is in $S^*$. Hence, we are unable to show that the sequence $\left(P_{\mathcal{C}}(y_k)\right)_{k\in \mathbb{N}}$ converges to the projected solution of QEP$(f,\mathcal{K})$ by applying Lemma \ref{lemma_Fejer_conv}. But, one can overcome this drawback by assuming the strict monotonicity of the bifunction $f$ as shown in the following result.
\begin{theorem}
Suppose $f: \mathcal{A} \times \mathcal{A} \to \mathbb{R}$ is strictly monotone and satisfies the assumptions \textbf{A} \ref{A1} and \textbf{A} \ref{A3}. Assume that $\mathcal{K}:\mathcal{C}\rightrightarrows \mathbb{R}^n$ satisfies assumption \textbf{A} \ref{A4} and  $S^{*}$ defined in \eqref{nonemptyset} is nonempty. Then the sequence $(y_k)_{k\in \mathbb{N}}$ produced by Algorithm \ref{algo} converges to $\bar y$ such that $P_{\mathcal{C}}(\bar y)$ is the projected solution of $QEP(f,\mathcal{K})$.   
\end{theorem}
\begin{proof}
    By proof of previous Theorem, we know that $ P_{\mathcal{C}}(\bar y)=\bar x\in S^P_{QEP} (f,\mathcal{K})$. Then, $\bar y\in \mathcal{K}(\bar x)$ satisfy
    \begin{equation}\label{haty}
        f(\bar y,y)\geq 0,\quad\forall\,y\in \mathcal{K}(\bar x).
    \end{equation} 
    We claim that $\tilde S_{QEP} (f,\mathcal{K})=S^*=\{\bar y\}$ where $\tilde S_{QEP} (f,\mathcal{K})$ is defined in \eqref{S_QEP}. % where $\bar x$ is the cluster point of the sequence $\{P_{\mathcal{C}}(y_k)\}$. 
    It is obvious that $S^*\subseteq \tilde S_{QEP} (f,\mathcal{K})$. We aim to show that $\tilde S_{QEP} (f,\mathcal{K})\subseteq S^*$.
    On contrary, let $\hat y\in S^*$ be arbitrary such that $\hat y\neq \bar y$. 
    Then, $\hat y\in \bigcap_{x\in \mathcal{C}} \mathcal{K}(x)$ satisfy  
     \begin{equation}\label{bary}
        f(\hat y,y)\geq 0,\quad\forall\,y\in \bigcup_{x\in \mathcal{C}} \mathcal{K}(x).
    \end{equation}
    Suppose $y=\hat y$ in (\ref{haty}), then we have $f(\bar y,\hat y)\geq 0.$ On the other hand, we obtain $f(\hat y,\bar y)\geq 0$ by taking $y=\bar y$ in (\ref{bary}). Hence, $f(\hat y,\bar y)+f(\bar y,\hat y)= 0$ due to monotonicity of $f$. But this contradicts the facts that $f$ is strictly monotone and $\hat y\neq \bar y$. Hence, our assumption that there is some $\hat y\in S^*$ such that $\hat y\neq \bar y$ is false. 
Now, the convergence of $(y_k)_{k\in \mathbb{N}}$ to $\bar y$ such that $P_{\mathcal{C}}(\bar y)$ is the projected solution follows from Lemma \ref{lemma_Fejer_conv} and Proposition \ref{Prop_Fejer_soln} (i).
\end{proof}

\section{Numerical Experiments}\label{num_exp}
 In this section, we present some computational results of Algorithms \ref{Algo3} and \ref{algo} on some examples. Both algorithms have been coded in MATLAB 2020a. The code has been executed on a personal computer with a processor INTEL Core i5 with 8GB RAM. We first solve an example with the help of Algorithm \ref{Algo3} and then the next two examples by using Algorithm \ref{algo}. 

\begin{figure}[H]
    \centering
    \includegraphics[scale=0.8]{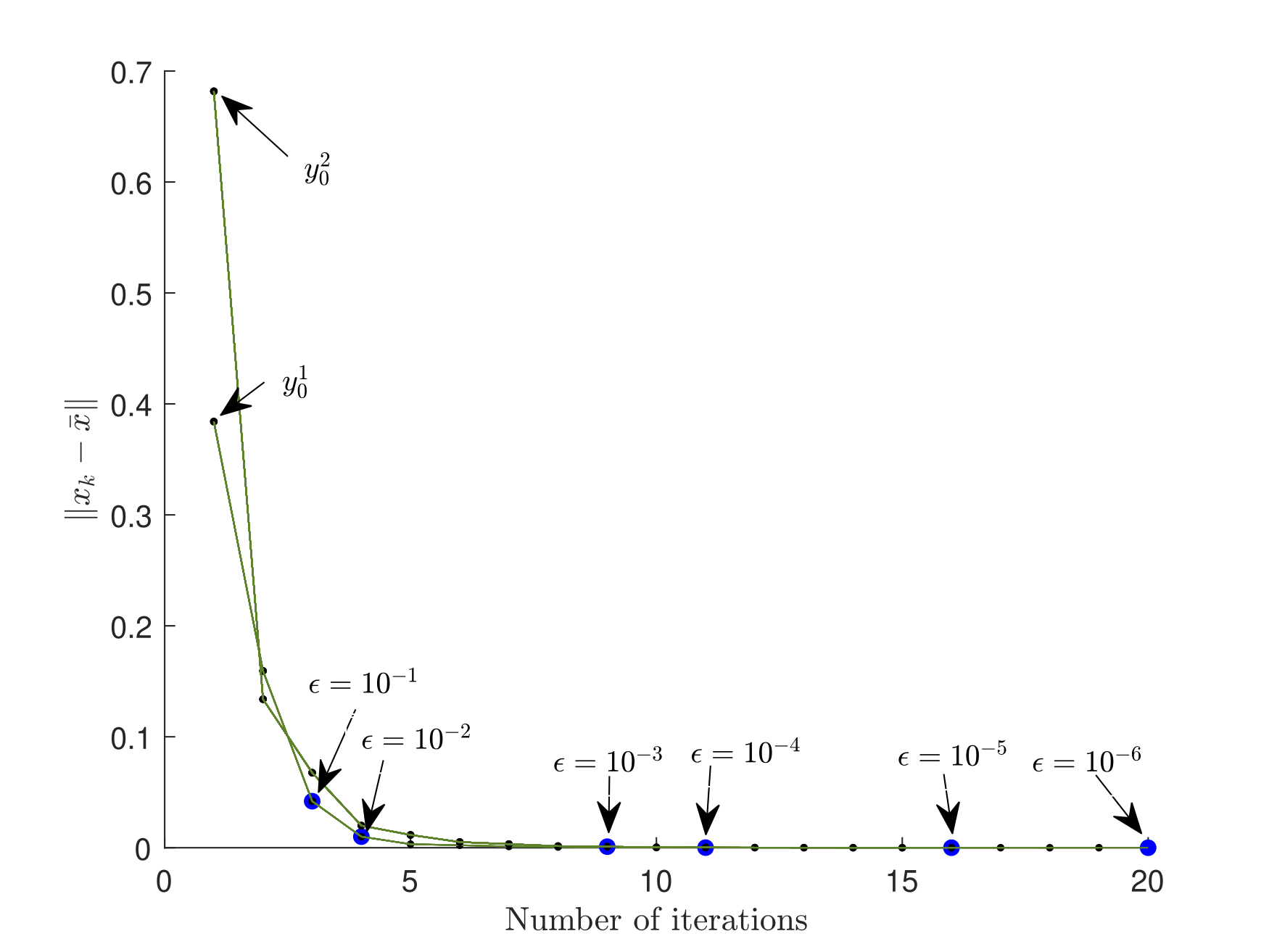}
    \caption{Comparison of convergence behavior of Algorithm \ref{Algo3} with number of iterations}
    \label{fig1_ex1}
\end{figure}

\begin{example}\label{ex_algo1}
Consider a bifunction $f:\mathbb{R}^5\times \mathbb{R}^5 \to \mathbb{R}$ of the form $f(x,y)=\langle F(x), y-x\rangle$, where $F(x) = \left(e^{-\|x\|^2}+0.1\right)\left(Mx+p\right)$ with $p=(-1,2,1,0,-1)^\top$ and
$$M=\begin{bmatrix}
5 & -1 & 2 & 0 & 2\\
-1 & 6 & -1 & 3 & 0\\
2 & -1 & 3& 0 & 1\\
0 & 3 & 0 & 5 & 0\\
2 & 0 & 1 & 0 & 4
\end{bmatrix}.$$
Let $ \mathcal{C} =\left\{z \in \mathbb{R}^{5}: \sum_{i=1}^{5}z_{i}\le 1, 0\le z_{i}\le 1\right\}$ and a set-valued map $\mathcal{K}:\mathcal{C} \rightrightarrows \mathbb{R}^{5}$ is defined as 
    $$\mathcal{K}(x) = -x +  \mathbb{B}\left[\bm{0},\frac{1}{2}\right],$$
    where $\mathbb{B}\left[\bm{0},\frac{1}{2}\right]$ is a closed ball with center origin and radius $\frac{1}{2}.$
\end{example}
It is worth noting that the bifunction $f$ in Example \ref{ex_algo1} is strongly monotone with $\mu = 0.1676$ and Lipschitz continuous with Lipschitz constant $L\approx 8.51$. Here, $\bar{x} = (0.096552,0,0,0.074467,0.109249)$ and $\bar{y} = (0.096552,-0.268553$, $-0.265827,0.074467,0.109249)$ are such that $P_{\mathcal{C}}(\bar{y}) = \bar{x}$ and $f(\bar{y},y)\ge 0$ for all $y \in \mathcal{K}(\bar{x}).$ Therefore, $\bar{x}$ is the projected solution to the above problem. The results obtained by Algorithm \ref{Algo3} with different initial points and accuracy level $\epsilon>0$ are presented in Table \ref{ex_1_table_alg_1} and Fig. \ref{fig1_ex1}. In Table \ref{ex_1_table_alg_1}, we also provide CPU time in seconds and total number of iterations (Iter) with respect to the initial point and accuracy level $\epsilon>0$ in columns 3 and 4, respectively.
%{p{4.3cm}|p{0.6cm}|p{1cm}|p{0.5cm}|p{4cm}}

\begin{table}[!h]
\centering
\caption{Obtained results by Algorithm \ref{Algo3} of Example \ref{ex_algo1} for different accuracy levels with different initial points}
\label{ex_1_table_alg_1}
\begin{tabular}{c|c|c|c|c}
\hline
\textbf{Initial point} & $\epsilon$ & \textbf{CPU time (in sec)} & \textbf{Iter} & \textbf{Projected solution} $(\bar{x})$ \\%& \textbf{Projected solution} $\bar{x}=P_{\mathcal{C}}(\bar{y})$ \\
\hline
 & $10^{-1}$  &  $25.39$  & $3$  & $(0.1327,0,0, 0.0539,0.1040)$ \\
  $y_{0}^{1}=(0.1,0.2,0.3,0,-0.1)$      & $10^{-2}$ & $38.86 $ & $4$ & $(0.1052,0,0, 0.0785,0.1122)$ \\
        & $10^{-4}$ & $115.05$ & $11$ & $(0.0965, 0,0, 0.0745,0.1092)$ \\   
\hline
 & $10^{-3}$  &  $86.95$  & $9$  & $(0.0970,0,0, 0.0749,0.1100)$ \\
  $y_{0}^{2}=(-2,0.6,-1.3,0.5,-1)$      & $10^{-5}$ & $165.16 $ & $16$ & $(0.0965,0,0, 0.0745,0.1092)$ \\
        & $10^{-6}$ & $209.61$ & $20$ & $(0.0966, 0,0,0.0745,0.1092)$ \\
\hline
\end{tabular}
\end{table}

In the next two examples, we test the performance of Algorithm \ref{algo} for different dimensions $n$ to find the projected solution.

\begin{example}\label{test_ex_1}
 Consider the nonempty, closed and convex set $\mathcal{C} = [-2, 0]^n$ and a set valued map $\mathcal{K}:\mathcal{C}\rightrightarrows \mathbb{R}^{n}$ is defined as 
\begin{equation*}
\mathcal{K} (x) = \prod_{i=1}^{n}\left[x_{i}-1,\frac{4x_{i}+11}{3}\right].
\end{equation*}
Let $f:\mathbb{R}^n\times \mathbb{R}^{n}\to \mathbb{R}$ be a bifunction defined as
$$f(x,y) = \sum_{i=1}^{n}\left(y_{i}^2-x_{i}^2\right) -2\sum_{i=1}^{n}\left(y_{i}-x_{i}\right). $$
\end{example}

\begin{figure}[H]
    \centering
    \includegraphics[scale=0.8]{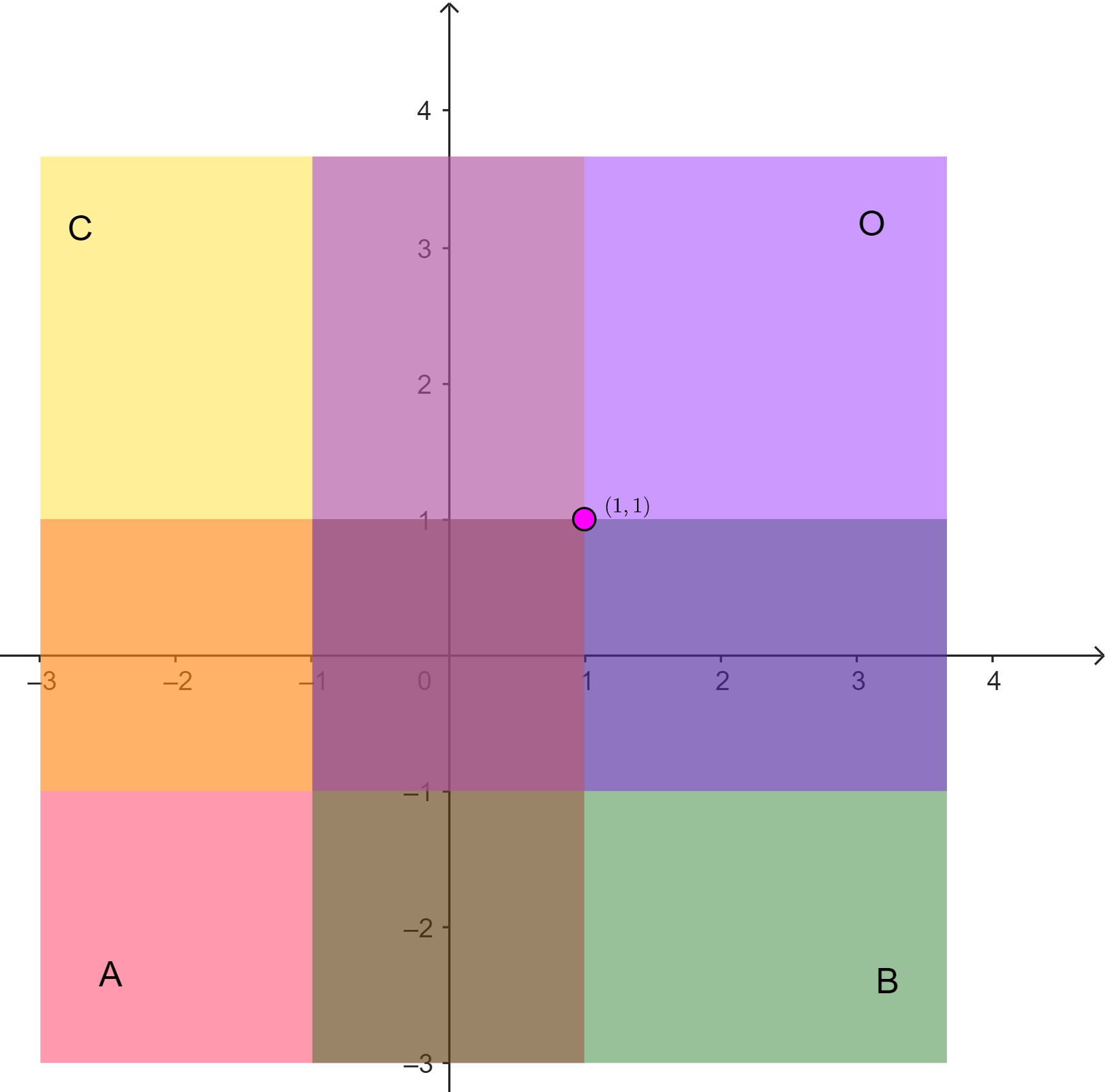}
    \caption{Illustration of the map $\mathcal{K}$ by taking four points from $\mathcal{C}$ for $n=2$  }
    \label{rectangles}
\end{figure}

When $n=2,$ then $\mathcal{C}=[-2,0]\times [-2,0]$ and $\mathcal{K}(\mathcal{C}) = [-3,\frac{11}{3}]\times [-3,\frac{11}{3}].$ The images of the points A$(-2,-2)$, B$(0,-2)$, C$(-2,0)$ and O$(0,0)$ of $\mathcal{C}$ are shown in Figure \ref{rectangles} in different colors. 

Note that the bifunction considered in Example \ref{test_ex_1} satisfies assumptions A\ref{A1}-A\ref{A4}. Furthermore,  $\bar{y}=(1,1\ldots,1)$ is such that $ \bar{y}\in \bigcap\limits_{z \in C}\mathcal{K}(z).$ Also, $f(\bar{y},y) \ge 0 $ for all $y \in \mathcal{K}(\mathcal{C}) = [-3, \frac{11}{3}]^n .$ Hence, $S^* = \{\bar{y}\}. $ In addition, $P_{\mathcal{C}}(S^*)=(0,0,\ldots,0),$ therefore, $S_{P}^{*} = \{(0,0,\ldots,0)\}$ is the subset of projected solution. To investigate the performance of Algorithm \ref{algo} to this example, we provide different starting points $(y_{0})$ and dimension $(n)$ to the algorithm and calculate the number of iterations (iter), CPU time with diferent sequence of $\gamma_{k}$ (see Table \ref{ex_1_table}).

\begin{table}[!h]
\centering
\caption{Performance of Algorithm \ref{algo} to Example \ref{test_ex_1} for different values of $n$ }
\label{ex_1_table}

\begin{tabular}{p{1cm}|p{2.5cm}|p{0.5cm}|p{1cm}|p{1.6cm}|p{1.6cm}|p{1.7cm}}
\hline
\textbf{$n$} & Initial point ($y_{0}$) & iter & CPU time (in sec) & $\gamma_{k}$&Obtained solution ($\bar{y}$) & Projected solution ($\bar{x}$) \\ \hline
$10$  &  $(2,2,\ldots,2)$ & $6$    & $0.0105$ &    $\frac{1}{k}$   &$(1,1,\ldots,1)$ &$(0,0,\ldots,0)$   \\
$100$  &  $(0.5,0.5,\ldots,0.5)$ & $3$    & $0.0542$ &    $\frac{1}{2k}$   &$(1,1,\ldots,1)$&$(0,0,\ldots,0)$ \\
$1000$  &  $(1.5,1.5,\ldots,1.5)$ & $5$    & $0.0255$ &    $e^{-k}$   &$(1,1,\ldots,1)$&$(0,0,\ldots,0)$ \\
$5000$  &  $(-1,-1,\ldots,-1)$ & $11$    & $0.0753$ &    $\frac{1}{k^2}$   &$(1,1,\ldots,1)$&$(0,0,\ldots,0)$ \\
$10000$  &  $(-5,-5,\ldots,-5)$ & $35$    & $0.3977$ &    1+$\frac{1}{\log(n+1)}$   &$(1,1,\ldots,1)$&$(0,0,\ldots,0)$ \\\hline
\end{tabular}
\end{table}

\begin{example}\label{test_ex_2}
Consider the quasi-equilibrium problem in which the set $\mathcal{C} = \left\{x\in \mathbb{R}^{n}_{+}: \sum\limits_{i=1}^{n} x_{i} \ge n\right\}$ and a multivalued mapping $\mathcal{K}$ is given by
$$\mathcal{K}(x) = \left\{y\in \mathbb{R}^{n}:\sum_{i=1}^{n}y_{i}\le 1+\frac{\|x\|}{1+\|x\|}, y_{i}\ge -\frac{1}{\|x\|}\right\}.$$
The bifunction $f:\mathbb{R}^n\times \mathbb{R}^{n}\to \mathbb{R}$ is of the form 
$$f(x,y) = \left(\sum_{i=1}^{n}y_{i}\right)^2-\left(\sum_{i=1}^{n}x_{i}\right)^2 - e^{\left(\sum\limits_{i=1}^{n}x_{i}-0.5\right)} + e^{\left(\sum\limits_{i=1}^{n}y_{i}-0.5\right)} +2\sum_{i=1}^{n}\left(x_{i}-y_{i}\right).$$
\end{example}
Clearly, the bifunction $f$ and the map $\mathcal{K}$ in Example \ref{test_ex_2} satisfies assumptions \textbf{A} \ref{A1}-\textbf{A} \ref{A4}. When $n=1,$ $\bar{y} = 0.5$ is such that $\bar{y}\in  \cap_{z\in \mathcal{C}}\mathcal{K}(z)=[0,1.5].$ Also, $f(\bar{y},y) \ge 0$ for all $y \in \cup_{z\in \mathcal{C}}\mathcal{K}(z)=[-1,2).$ Therefore, $S^* \neq \emptyset$ for $n=1$ and hence $S_{P}^{*} \neq \emptyset.$  For $n=2,$ $\cap_{z\in \mathcal{C}}\mathcal{K}(z)=\{y\in\mathbb{R}^2: 0<y_1+y_2\le 3-\sqrt{2}\}$ and $\cup_{z\in \mathcal{C}}\mathcal{K}(z)= \{y\in \mathbb{R}^2:-\sqrt{2}\le y_1+y_2<2\}.$ Also, $\bar{y} = (0.3675,0.3675) \in \cap_{z\in \mathcal{C}}\mathcal{K}(z)$ and $f(\bar{y},y) \ge 0$ for all $y \in \cup_{z\in \mathcal{C}}\mathcal{K}(z).$ The numerical results of this example obtained by Algorithm \ref{algo} are presented in Table \ref{ex_31_table} taking different values of $n$ and initial points. It will not be possible to verify that for each $n$, the obtained projected solutions belong to the set $S_{P}^*$. But, we have performed the proposed algorithm for some values of $n(>2)$ and the obtained results are added in Table \ref{ex_31_table}.

% \begin{table}[H]
% \centering
% \caption{} 
% \begin{tabular}{c c c c c c}
% \toprule
% \textbf{$n$} & $C$ & $\cap_{z\in C}K(z)$ & $\cup_{z\in C}K(z)$ & $f(x,y)$& $\bar{y}\in S^*$ \\ \midrule
% $1$  &  $[1,\infty)$ & $[0,1.5]$    & $[-1, 2)$ &    $f$   &$0.5$   \\
% $2$  &  $(0.5,0.5)$ & $3$    & $0.0844$ &    $\frac{1}{2k}$   &$(0.3675,0.3675)$  \\
% $100$ & $(0.5,\ldots,0.5)$ &$12$ & $ 4.1913$ & $\frac{1}{5k}$ & $(0.0132,\ldots,0.0132)$\\
% \bottomrule
% \end{tabular}
% \end{table}

\begin{table}[!h]
\centering
\caption{Performance of Algorithm \ref{algo} to Example \ref{test_ex_2} for different values of $n$ }
\label{ex_31_table}
\begin{tabular}{p{0.4cm}| p{1.8cm}| p{0.4cm}| p{1cm}| p{0.7cm}| p{2.9cm}| p{2.0cm}}
\hline
\textbf{$n$} & Initial point ($y_{0}$) & iter & CPU time (in sec) & $\gamma_{k}$&Obtained solution ($\bar{y}$) & Projected solution ($\bar{x}$) \\ \hline
$1$  &  $2$ & $5$    & $0.1491$ &    $\frac{1}{k}$   &$0.5$ & $1$  \\
$2$  &  $(0.5,0.5)$ & $3$    & $0.0844$ &    $\frac{1}{2k}$   &$(0.3675,0.3675)$ & $(1,1)$  \\
$5$ &  $(1,\ldots,1)$ &$3$ &$0.1015$ & $\frac{1}{3k}$  &$(0.1952,\ldots,0.1952)$& $(1,\ldots,1)$ \\
$10$ & $(2,\ldots,2)$ & $3$ &$0.1924$ & $\frac{1}{4k}$ & $(0.1079,\ldots,0.1079)$ & $(1,\ldots,1)$ \\
$100$ & $(0.5,\ldots,0.5)$ &$12$ & $ 4.1913$ & $\frac{1}{5k}$ & $(0.0132,\ldots,0.0132)$ & $(1,\ldots,1)$\\ \hline
\end{tabular}
\end{table}

\section{Application to Generalized Nash Equilibrium Problems}\label{sec_GNEP}
Arrow-Debreu \cite{debreu} extended the well-known Nash equilibrium problems (NEP) to the case where the strategy sets of players are allowed to rely on the strategies selected by rivals. These problems are known as generalized Nash equilibrium problems (GNEP) (or ``abstract economy"). Suppose any player $i$ in $\Lambda=\{1,2,\cdots N\}$ controls a strategy vector $x_i\in C_i\subseteq \mathbb{R}^{n_i}$ where $\sum_{i\in \Lambda} n_i=n$. Assume that,
				\begin{equation*}\label{X}
					C_{-i}=\prod_{k\neq i} C_k,\quad \mathcal{C}=\prod_i C_i,
				\end{equation*}
	such that the vector $x_{-i}\in C_{-i}\subseteq \mathbb{R}^{n-n_i}$ is strategy vector for all players except $i$ and $x=(x_i)_{i\in \Lambda}=(x_i,x_{-i})\in \mathcal{C}$ indicates the strategies of all $N$-players. Suppose $\theta_i:\mathbb{R}^n\rightarrow \mathbb{R}$ is the objective/profit function of player $i$. For the given strategy vector $x_{-i}$ of rivals, the player $i$ intends to choose a strategy $x_i$ in the feasible strategy set $K_i(x)\subseteq \mathbb{R}^{n_i}$ such that $x_i$ maximizes the function $\theta_i(\cdot,x_{-i})$. Then, a vector $\bar x\in \mathcal{C}$ is known as generalized Nash equilibrium if			\begin{equation}\label{GNEP2}
			\bar x_i\in K_i(\bar x)~\text{and}~	\bar x_i\in \arg \max_{K_i(\bar x)} \theta_i(\cdot,\bar x_{-i}),\quad\forall\, i\in \Lambda.
			\end{equation}
            The problem of finding $\bar x\in \mathcal{C}$ satisfying \eqref{GNEP2} is known as generalized Nash equilibrium problem denoted by GNEP$(\theta_i,K_i)_{i\in \Lambda}$. If $K_i(x)=D_i$ for some, $D_i\subseteq C_i$ for all $i\in \Lambda$ then GNEP$(\theta_i,K_i)_{i\in \Lambda}$ reduces to Nash equilibrium problem NEP$(\theta_i,D_i)_{i\in \Lambda}$. Note that most of the existence result on GNEP$(\theta_i,K_i)_{i\in \Lambda}$ in literature assumes that the constraint map $\mathcal{K}: \mathcal{C}\rightrightarrows \mathbb{R}^n$ defined as $\mathcal{K}(x)=\prod_i K_i (x)$ is a self map, i.e., $K_i(\mathcal{C})\subseteq C_i$ for all $i\in \Lambda$ (see e.g., \cite{debreu,facc_survey}). However, this need not be true always as illustrated in \cite{aussel2016existence} with the help of electricity market model. Let us recall the concept of projected solutions from \cite{aussel2016existence} for a GNEP$(\theta_i,K_i)_{i\in \Lambda}$ in which $\mathcal{K}$ need not be a self-map. A vector $\bar x\in \mathcal{C}$ is known as projected solution for GNEP$(\theta_i,K_i)_{i\in \Lambda}$ if there exists $\bar{y}\in \mathcal{K}(\bar{x})$ such that
\begin{itemize}
\item[(a)] $\norm{\bar{x}-\bar{y}}=\inf_{y\in \mathcal{C}}  \norm{\bar{x}-{y}}$;
\item[(b)] $\bar y$ solves NEP$(\theta_i,K_i(\bar x))_{i\in \Lambda}$.
	\end{itemize}  

    It is well-known that the set of projected solutions for GNEP$(\theta_i,K_i)_{i\in \Lambda}$ coincides with the set of projected solutions for an auxiliary quasi-variational inequality problem. \begin{lemma} \cite[Lemma 4.1]{aussel2016existence}\label{QVI_reform_aussel}
				Consider a GNEP$\,(\theta_i, K_i)_{i\in \Lambda}$ in which $C_i$ is a non-empty closed convex set and $K_i(x)$ is non-empty closed convex for any $x\in \mathcal{C}$. Let $\theta_i$ be continuously differentiable and the function $\theta_i(\cdot,x_{-i})$ be concave for any $x_{-i} \in \mathbb{R}^{n-n_i}$. Then $\bar x\in \mathcal{C}$ is a projected solution for $GNEP(\theta_i,K_i)_{i\in \Lambda}$ if and only if it is projected solution for QVI$(G,\mathcal{K})$ where $G(x)=(-\nabla_{x_i} \theta_i(x))_{i\in \Lambda}$ and $\mathcal{K}(x)=\prod_i K_i(x_{-i})$.
				\end{lemma}
                
    Suppose $f:\mathbb{R}^n\times \mathbb{R}^n\rightarrow \mathbb{R}$ is defined as 

\begin{equation}\label{define_f_GNEP}
    f(y,z)=\sum_{i\in I}\langle -\nabla_{y_i} \theta_{i}(y),z_i-y_i\rangle
\end{equation} and $\mathcal{K}:X\rightrightarrows \mathbb{R}^n$ is defined as $\mathcal{K}=\prod_{i\in I} K_i$. It is easy to observe from Lemma \ref{QVI_reform_aussel} that any $\bar x$ is a projected solution of QEP($f,\mathcal{K}$) iff it is a projected solution of GNEP$(\theta_i,K_i)_{i\in \Lambda}$. Hence, we ensure the convergence of the sequence generated by Algorithm \ref{Algo3} to a projected solution of GNEP$(\theta_i,K_i)_{i\in \Lambda}$ by considering the set $\mathcal{C}=\prod_{i\in \Lambda} C_i$, the map $\mathcal{K}=\prod_{i\in \Lambda} K_i$, the set $\mathcal{A}=\overline{co}(\mathcal{K}(\mathcal{C}))$ and the bifunction $f$ as \eqref{define_f_GNEP} in Theorem \ref{Thm_strong_mono}.
\begin{corollary}\label{Coro_GNEP_Str_mono}
Consider a GNEP$\,(\theta_i, K_i)_{i\in \Lambda}$ in which $C_i$ is a non-empty closed convex set and  $K_i(x)$ is non-empty closed convex for any $x\in \mathcal{C}$. %Suppose $f:\mathbb{R}^n\times \mathbb{R}^n\rightarrow \mathbb{R}$ is defined as (\ref{define_f_GNEP}) and $\mathcal{K}=\prod_{i\in I} K_i$ where $K_i:X\rightrightarrows \mathbb{R}^{n_i}$ admits non-empty closed convex values. 
Assume that:
    \begin{itemize}
        \item [(i)] $\theta_i$ is continuously differentiable function and $\theta_i(\cdot,x_{-i})$ is concave for any $x_{-i}$;
        \item[(ii)] there exists $L>0$ such that $\norm{G(y)-G(z)}\leq L \norm{y-z}$ for all  $y,z\in \mathcal{A}$ where $G(x)=(-\nabla_{x_i} \theta_i(x))_{i\in \Lambda}$;
        \item [(iii)] $f$ is strongly monotone bifunction with parameter $\mu>0$;
        \item[(iv)] there exists $\alpha\in (0,\frac{\mu}{L})$ such that $\|P_{\mathcal{K}(u)} (z)-$$P_{\mathcal{K}(v)}(z)\|$$\leq \alpha\|u-v\|$ for all $u,v\in \mathcal{C}$ and $z\in \mathbb{R}^n$.       
    \end{itemize}
    Then, there exists a unique projected solution $\bar x$ for GNEP$(\theta_i,K_i)_{i\in \Lambda}$. Further, the sequence generated by Algorithm \ref{Algo3} converges to this unique projected solution $\bar x$ in the following way:
    \begin{equation*}
        \norm{x_k-\bar x}\leq \frac{1}{\delta} \exp \left \{\frac{-\delta}{2}k\right\} \norm{y_0-T(x_0)},
    \end{equation*}
    where $\delta=1-\alpha\frac{L}{\mu}$.
\end{corollary}
\begin{remark}
Note that hypothesis (iii) of the above result is easy to check if $\theta_i$ is twice differentiable. In fact, the strong monotonicity parameter $\mu$ of $f$ is equal to the smallest eigenvalue of the Jacobian of $F(y)=(-\nabla_{y_i} \theta_i(y))_{i\in \Lambda}$ (see \cite[Theorem 5.4.3]{book}).
\end{remark}
\begin{remark}
The hypothesis (ii) of the above result is easy to verify if $\theta_i$ is twice differentiable. Suppose $J_F(x)$ denotes the Jacobian of $F(x)$ at any $x\in \mathcal{C}$. Then the Lipschitz continuity constant $L$ of $F$ is given by $\norm{J_F}$. Suppose $\lambda(x)$ denotes the largest eigen value of $(J_F(x))^T J_F(x)$ for any $x\in \mathcal{C}$. Then $\norm{J_F}$ is equal to the positive square root of the supremum of $\lambda(x)$ over $\mathcal{C}$ (Lemma 3.1 in \cite{Khalil}).
\end{remark}
We show the convergence of the sequence generated by Algorithm \ref{algo} to the projected solution of GNEP$\,(\theta_i,\mathcal{K}_i)_{i\in \Lambda}$ by considering the set $\mathcal{C}=\prod_{i\in \Lambda} C_i$, the map $\mathcal{K}=\prod_{i\in \Lambda} K_i$, the set $\mathcal{A}=\overline{co}(\mathcal{K}(\mathcal{C}))$ and the bifunction $f$ as \eqref{define_f_GNEP} in Theorem \ref{Thm_proximal_conv}. %Suppose $f:\mathbb{R}^n\times \mathbb{R}^n\rightarrow \mathbb{R}$ is defined as (\ref{define_f_GNEP}) and $\mathcal{K}=\prod_{i\in I} K_i$ where $K_i:X\rightrightarrows \mathbb{R}^{n_i}$.
\begin{corollary} \label{corollary_GNEP_PPA}
Consider a GNEP$\,(\theta_i, K_i)_{i\in \Lambda}$ in which $C_i$ is a non-empty closed convex set. Assume that:
    \begin{itemize}
        \item [(i)] $\theta_i$ is continuously differentiable function and $\theta_i(\cdot,x_{-i})$ is concave for any $x_{-i}$;
        \item [(ii)] $K_i$ is closed lower semi-continuous map such that $K_i(x)$ is non-empty closed convex for any $x\in \mathcal{C}$;
        \item[(iii)] $f$ is monotone bifunction such that for any sequence $(y_n)_{n \in \mathbb{N}}\subset \mathcal{A}$ with $\lim_n\norm{y_n}=\infty$, there exists $z\in \mathcal{A}$ and $n_0\in \mathbb{N}$ such that $f(y_n,z)\leq 0$ for all $n\geq n_0$;
        \item[(iv)] the set $S^*$ defined in (\ref{nonemptyset}) is non-empty.      
    \end{itemize}
    Then, every cluster point $\bar y$ of the sequence $(y_k)_{k\in \mathbb{N}}$ generated by Algorithm \ref{algo} yields the projected solution $\bar x=P_{\mathcal{C}}(\bar y)$ of GNEP$\,(\theta_i,K_i)_{i\in \Lambda}$.
\end{corollary}

\section{Numerical approach to Electricity Market Model}\label{EMM_NUM}
In this section, we aim to find the projected solutions for the deregulated electricity market model (EMM) proposed in \cite{aussel2016existence}. Let us recall the electricity market model from \cite[Section 2.2]{aussel2016existence}. By following the electricity market model presented in \cite{aussel2016existence}, we assume there are $N$ producers in the set $\Lambda=\{1,\cdots,N\}$ acting as leaders and an independent system operators (ISO) acting as follower in a multi-leader-common-follower game (MLFCG). Each producer wants to maximize profits, and ISO wants to minimize the overall price of electricity. Assume that the demand $D$ is given in advance. Suppose $Q_i$ indicates the production capability of any producer $i$ and $q_i\in [0,Q_i]$ indicates his production quantity. Each producer proposes a unit price bid function $\psi_i:[0,Q_i]\rightarrow \mathbb{R}$. %and, which is directly proportional to quantity produced \cite{hawthorne}. 
For $q_i$ units, suppose $\phi_i(q_i)=\int_{0}^{q_i} \psi_i(q) dq+r_i$ and $P_i(q_i)$ evaluate the revenue and production cost of $q_i$ units, respectively. For given bid functions $\psi_i, i\in \Lambda$, the aim of ISO is to determine $q_i$ so that the overall cost charged by producers is minimized and demand is also fulfilled. Each producer aims to find a pair of bid function and production quantity $(\psi_i,q_i)$ in such a way that profit $\phi_i(q_i)-P_i(q_i)$ is maximized.

Clearly, each producer $i\in \Lambda$ intends to solve the following problem:
\begin{align*}
	&\max_{\phi_i,q_i}~ [\phi_i(q_i)-P_i(q_i)]\\
	&~\text{s.t.}~ \phi_i~\text{is admissible bid function satisfying}~\text{ISO}(\phi);
\end{align*}
and the problem $\text{ISO}(\phi)$ is given as:
\begin{align*}
	&\min_q \sum_{i\in \Lambda}\phi_i(q_i)\\
	&~ \text{s.t.}~q_i\in [0,Q_i]~\text{satisfies} ~\sum_{i\in \Lambda} q_i\geq D.
\end{align*}

According to the existing literature, the production cost $P_i$ is a quadratic function of $q_i$ (see \cite{aussel2016existence,ausselnash1,ausselnash2}).  Suppose $P_i(q_i)=A_i q_i^2+B_i q_i$ the parameters $A_i>0$ and $B_i\in \mathbb{R}$ are fixed and known to producers. Suppose $0=q_i^0<q_i^2<\cdots< q_i^{k}=Q_i$ forms the partition of $[0,Q_i]$ for a fixed $k\geq 2$. Then, the unit price bid function $\psi_i$ is an increasing step function as follows:
	\begin{align*}\label{psi}
		\psi_i(q_i)=\begin{cases} 
			p_i^0,&\text{if}~q_i=0\\
			p_i^{m},&\text{if}~q_i\in (q_i^m,q_i^{m+1}],~\text{for all}~m=0,1,\cdots, k-1
		\end{cases}
	\end{align*}
	where $p_i^0<p_i^1\cdots<p_i^{k-1}$. For any pair of producers $i,j\in \Lambda$, let us consider $q_i^m=q_j^m$ for all $m\in \{0,1,\cdots,k\}$. Then, the resulting piece-wise linear revenue bid function $\phi_i: [0,Q]\rightarrow \mathbb{R}$ is defined as,
		\begin{equation*}
			\phi_i(q)=\begin{cases}
				p_i^0 q+p_i^0,~&\textit{if}~q\in[q_0,q_1],\\
				p_i^1q+p_i^0+(p_i^0-p_i^1)q_1,~&\textit{if}~q\in[q_1,q_2],\\
				p_i^2 q+p_i^0+(p_i^0-p_i^1)q_1+(p_i^1-p_i^2)q_2~&\textit{if}~q\in[q_2,q_3],\\\vdots\\
				p_i^{k-1}q+p_i^0+\sum_{n=1}^{k-1}(p_i^{n-1}-p_i^n)q_n~&\textit{if}~q\in[q_{k-1},q_k].
				\end{cases}
		\end{equation*}

        For any $j=0,\cdots,k-1$, suppose $p_{i,j}^{min}$ and $p_{i,j}^{max}$ are fixed lower and upper bounds of $p_i^j$ with $p_{i,j}^{max}<p_{i,j+1}^{min}$. Then the set 
        \begin{small}
        \begin{align*}
			C_i=\{\phi_i\in L^2([0,Q],\mathbb{R})\,|&~\phi_i(q)=\alpha_{j} q+\beta_{j},q\in [q_{j},q_{j+1}]~\text{for}~ j=0,\cdots,k-1,~\text{where}~\alpha_{j}=p_i^{j},\\&~p_i^j\in [p_{i,j}^{min},p_{i,j}^{max}],\beta_{j-1}=p_i^0+\sum_{n=1}^{j}(p_i^{n-1}-p_i^n)q_n\}
	\end{align*} 
    \end{small}
        is closed and convex. Further, any $\phi\in C_i$ satisfies $p_i^0<p_i^1\cdots<p_i^{k-1}$
	 We know that $C_i$ is finite dimensional subset of $L^2([0,Q],\mathbb{R})$ isometric to subset of $\mathbb{R}^{3k}$. Suppose $\mathcal{C}=\prod_{i=1}^N C_i$. Let us define $I:C_i\rightarrow \mathbb{R}^{3k}$ as,
	\begin{equation*}
		I(\phi_i)=(0,\cdots,0, \alpha_i^0,\alpha_i^1,\cdots,\alpha_i^{k-1},\beta_i^0,\beta_i^1,\cdots,\beta_i^{k-1}).
	\end{equation*}
    Suppose the quadratic bid function of any producer $y_i= a_iq_i^2+b_iq_i+c_i$ satisfies $a_i=A_i$, $c_i\geq p_i^1$, $c_i\in[\underbar c_i,\bar c_i]$ and $b_i\in [\underbar b_i,\bar b_i]$ where $\underbar b_i,\bar b_i,\underbar c_i,\bar c_i$ are fixed non-negative quantities and $\bar c_i\geq p_{i,0}^{max}$. Then the map $K_i:C\rightrightarrows L^2([0,Q],\mathbb{R})$ is,
	\begin{align}\label{define_map_ki}
	K_i(\phi)=\{y_i\,|\,y_i(q)=a_i q_i^2+b_i q_i+c_i, a_i=A_i,b_i\in [\underbar b_i,\bar b_i],c_i\geq p_i^0, c_i\in[\underbar c_i,\bar c_i]\}.
	\end{align}
	 Then $K_i(\mathcal{C})$ is finite dimensional subset of $L^2([0,Q],\mathbb{R})$ isometric to subset of $\mathbb{R}^{3k}$.
     Let us define $I:K_i(\mathcal{C})\rightarrow \mathbb{R}^{3k}$ as,
		\begin{align*}
			I(y_i)=(a_i,\cdots,a_i,b_i,\cdots,b_i,c_i,\cdots,c_i).
		\end{align*}
	Since $I$ is isometry, for any $y_i\in K_i(\mathcal{C})$ and $\phi_i\in C_i$, we can write
		\begin{equation}\label{distance}
			d(I(y_i),I(\phi_i))_{\mathbb{R}^{3k}}=\norm{y_i-\phi_i}_{L^2}.
		\end{equation}
Hence, the problem of finding a projected solution for a given EMM can be rewritten in $\mathbb{R}^{3Nk}$.
    % Suppose the quadratic bid function of any producer $y_i= a_iq_i^2+b_iq_i+c_i$ satisfies $a_i=A_i$, $c_i\geq p_i^1$ and $b_i\in [v_i^1,v_i^2]$ where $v_i^1$ and $v_i^2$ are fixed non-negative quantities.
    Under the assumption that the production quantity $q_i$ by any producer $i$ satisfies $q_i\in (0,Q_i)$, Aussel et. al. \cite{aussel2016existence} observed that the given electricity market is a particular case of GNEP$(\theta_i,K_i)_{i\in \Lambda}$ where $\theta_i(b)=(b_i-B_i) q_i(b)+c_i$ for
\begin{equation}\label{exp_qib}
    q_i(b)=\frac{1}{2A_{i}}\left[\frac{2\left(\prod_{k=1}^{10}A_k\right) + \sum_{j=1}^{10}\left(\prod_{k\neq j}A_{k}\right)b_{j}}{\sum_{j=1}^{10}\left(\prod_{k\neq j}A_k\right)}-b_i\right].
\end{equation}
and $K_i$ is given by \eqref{define_map_ki}.
Our aim is to find a vector of bid functions $\tilde\phi=(\tilde \phi_i)_{i=1}^n$ in $\mathcal{C}=\prod_{i=1}^{n} C_i$ such that there exists a vector $\tilde y=(\tilde y_i)_{i\in \Lambda}$ formed by quadratic functions $\tilde y_i(q)=A_i q^2+b_iq+c_i$ (also expressed as $I(\tilde y_i)=(\tilde A_i,\tilde b_i,\tilde c_i)$) which satisfies:
		\begin{itemize}
			\item[(a)] $\tilde \phi=(\tilde \phi_i)_{i\in \Lambda}$ is projection of $\tilde y=(\tilde y_i)_{i\in \Lambda}$ on $\mathcal{C}$, that is, $\tilde \phi$ satisfies
\begin{align}\label{projectedprob1}
			\inf_{\phi \in \mathcal{C}} d(I(\tilde y),I(\phi))_{\mathbb{R}^{3Nk}};%\sum_{i\in \Lambda} \int_{0}^{Q_i} |\tilde y_i(q_i)-\phi_i(q_i)|^2 dq_i;
			\end{align}
            where the metric $d$ is defined as \eqref{distance}.
			\item[(b)]  Assuming that $b=(\tilde b_{-i},b_i)$, the vector $(\tilde b_i,\tilde c_i)$ solves the following problem for each $i\in\Lambda$,			
			\begin{equation}\left. 
			\begin{aligned} 
			&\max_{(b_i,c_i)}~ [(b_i-B_i)q_i(b)+c_i]\\
			&~\text{s.t.}~I(y_i)=(A_i,b_i,c_i)\in I(K_i(\tilde \phi))~\text{and}~ q_i(b)~\text{is given by}~(\ref{exp_qib}).\label{projectedprob2}
			\end{aligned}\right\}
			\end{equation}
		\end{itemize}  
        The vector $(\tilde \phi,\tilde y)$ satisfying \eqref{projectedprob1} and \eqref{projectedprob2} are known as projected solution pair for EMM.
\subsection{For 2 Producers by using Algorithm \ref{algo}}
For $2$ producers, suppose $y_i=(A_i,\cdots,A_i,b_i,\cdots,b_i,c_i,\cdots,c_i)\in \mathbb{R}^{3k}$ and $D=1$. Then, from \cite{aussel2016existence} we have 
	\begin{align*}
		\theta_1(y_1,y_2)=(b_1-B_1)q_1(b)+c_1,\theta_2(y_1,y_2)=(b_2-B_2)q_2(b)+c_2,
	\end{align*}
	where,
	$q_i(b)=\frac{1}{2A_i}\bigg[\frac{2A_1 A_2+A_{-i}b_i+A_{i}b_{-i}}{A_1+A_2}-b_i\bigg]$.
    
    Suppose $k=2$ and $q_0=0,q_1=\frac{1}{2},q_2=1$. Let $\mathcal{A}=I(K(C))\subset \mathbb{R}^{12}$ and $y,z\in \mathcal{A}$ has the form $y=(A_i,\cdots,A_i,b_i,\cdots,b_i,c_i,\cdots,c_i)_{i=1}^2$ and $z=(A_i,\cdots,A_i,v_i,\cdots,v_i,w_i,\cdots,w_i)_{i=1}^2$. Then, we form the bifunction $f:\mathcal{A}\times \mathcal{A}\to \mathbb{R}$ as $f(y,z)= \sum_{i=1}^2\langle\nabla_{y_i}\theta_i(y),y_i-z_i\rangle$, that is,
    \begin{align}
        f(y,z)=&\frac{1}{2(A_1+A_2)}\left[(b_1-v_1)(2A_1+B_1-2b_1+b_2)+(b_2-v_2)(2A_1+B_2+b_1-2b_2)\right]\notag\\&+(c_1+c_2-w_1-w_2).\label{bifunction_appl}
    \end{align}
\begin{proposition}\label{Prop_appl_mono}
    The bifunction $f$ defined as (\ref{bifunction_appl}) is monotone.
\end{proposition}
\begin{proof}
   Let $G(y)=(-\nabla_i\theta_i(y))_{i=1}^2$. One can observe that there are 3 different eigenvalues of the Jacobian of $G(y)$ that can be given as $0,0,\frac{1}{2(A_1+A_2)}, \frac{3}{2(A_1+A_2)}$. Hence $G:\mathcal{A}\rightarrow \mathbb{R}^{12}$ is monotone due to \cite[Theorem 5.4.3]{book}. Then, it is easy to check that the bifunction $f$ is also monotone.
\end{proof}
Following result derives the conditions under which $S_P^*$ is non-empty.
\begin{proposition}\label{prop_A,b,c}
    If $A_i, B_i,\bar b_i,\bar c_i$ are positive quantities for $i=1,2$ and anyone of the following conditions holds then $S_P^*\neq \emptyset$,
    \begin{itemize}
        \item [(i)] $2\bar b_1\leq \bar b_2$ and $2A_1+B_2+\bar b_1-2\bar b_2\geq 0$;
        \item [(ii)] $\bar b_1\geq 2 \bar b_2$ and $2A_2+B_1-2\bar b_1+\bar b_2\geq 0$.
    \end{itemize}
\end{proposition}
\begin{proof}
   % We know that $ f(y,z)=\frac{1}{2(A_1+A_2)}\bigg[(b_1-v_1)(2A_1+B_1-2b_1+b_2)$$+(b_2-v_2)(2A_1+B_2+b_1-2b_2)\bigg]$.
   In the case condition (i) holds, we have $f(\bar y,z)\geq 0$ for any $z\in \mathcal{K}(\mathcal{C})$ by taking $\bar y=(A_i,\bar b_i,\bar c_i)_{i=1}^2$. Further $\bar y\in \cap_{\phi\in {\mathcal C}} \mathcal{K}(\phi)$. Hence, the set $S^*\neq \emptyset$, which further implies that the set $S_P^*\neq \emptyset$.
\end{proof}
The following result is a consequence of Corollary \ref{corollary_GNEP_PPA} in the view of Proposition \ref{Prop_appl_mono}.
\begin{proposition}\label{Prop_Appl_EMM}
    Suppose any one of the conditions (i) or (ii) in Proposition \ref{prop_A,b,c} hold. Then, any cluster point $\tilde y$ of the sequence $(y_k)_{k\in \mathbb{N}}$ generated by Algorithm \ref{algo} yeilds a projected solution $\tilde \phi=P_{\mathcal{C}} (\tilde y)$ of the given electricity market.
\end{proposition}
We have performed Algorithm \ref{algo} on EMM with two producers by taking $p_1^0\in [15,50]$, $p_2^0\in [20,60]$, $p_1^1\in [51,200]$ and $p_2^1\in[61,100]$. The results of the EMM obtained by applying Algorithm \ref{algo} are presented in Table \ref{EMM_table} with accuracy level $\varepsilon = 10^{-4}$. In the first two instances of Table \ref{EMM_table}, the conditions of Proposition \ref{Prop_Appl_EMM} are true. Further, in the instance 3, the conditions of Proposition \ref{Prop_Appl_EMM} does not hold.

\begin{table}[!h]
\centering
\caption{Performance of Algorithm \ref{algo} to electricity market model for two producers}
\label{EMM_table}
\begin{tabular}
{p{1.4cm}|p{1.8cm}|p{2.5cm}|p{0.4cm}|p{0.9cm}|p{2.8cm}|p{1.6cm}}
\hline
 Coeffs of cost function & Bounds of $b_i$ and $c_i$  & Initial point $y_{0}=(A_i^0,b_i^0,c_i^0)_{i=1}^2$ & iter & CPU time (in sec) &  $(p_{i}^{0},p_{i}^{1})_{i=1}^2$ & Conclusion \\ \hline
$A_1$=10, $A_2=5,$& $b_1\in$[0,10], $b_2\in$[0,20],  & (10,10,16,5,49,46) & $40$    & $69.86$  & \multirow{3}{4em}{(37.26,51,45.97,61)}  & \multirow{4}{4em}{(Projected solution)} \\
$B_1$=2, $B_2=10$ & $c_1\in$[15,50], $c_2\in$[20,60] & (10,0,17,5,48,47) & $41$    & $109.18$  &   \\ \hline
 $A_1$=6, $A_2$=18, &$b_1\in$[0,31], $b_2\in$[0,15] &  (6,5,11,18,17,42) & $20$    & $187.02$    & \multirow{3}{4em}{(50,173.75,60,78.12)} & \multirow{4}{4em}{(Projected solution)}    \\
$B_1$=12, $B_2=5$&$c_1\in$[15,100], $c_2\in$[20,100]  &  (6,0,0,18,0,0) & $25$ & $244.38$  &  &\\ \hline
 $A_1$=6, $A_2=7,$ & $b_1\in$[0,31], $b_2\in[0,50]$& (6,0,0,7,0,0)& $69$ & $335.60$ &\multirow{3}{4em}{(50,173.75,60,100)} & \multirow{4}{4em}{(No conclusion)} \\
$B_1$=100, $B_2=20$ & $c_1\in$[15,100], $c_2\in$[20,100] & (6,1,10,7,2,5) &$68$&$556.39$&   & \\
\hline
\end{tabular}
\end{table}

\subsection{For 10 Producers by using Algorithm \ref{Algo3}}
In the previous subsection, the projected solution for the electricity market model is obtained by employing Proposition \ref{prop_A,b,c}. For more than $2$ producers, the condition derived in Proposition \ref{prop_A,b,c} is not easy to verify. Hence, we have employed Algorithm \ref{Algo3} to find the projected solution in case of $10$ producers. 

A constraint map $\mathcal{K}$ satisfying (\ref{map_K}) is also known as moving set constraint map. This type of map fulfills hypothesis \ref{hy_3.1ii} of Theorem \ref{Thm_strong_mono}. In the above electricity market model, if we have $I(\mathcal{K}):\mathcal{C}\rightrightarrows \mathbb{R}^{30k}$ as,
\begin{equation*}\label{define_map_Ki}
	I(\mathcal{K}(\phi))=\{I(y)=(A_i,b_i,c_i)_{i=1}^{10}|y_i(q)=A_i q_i^2+b_i q_i+c_i,b_i\in [\underbar b_i,\bar b_i],c_i=\alpha p_i^0,\,\text{for all}\,i\}.
\end{equation*}
Then, $\mathcal{K}$ is the moving set constraint map that satisfies the hypothesis \ref{hy_3.1ii} of Theorem \ref{Thm_strong_mono} by choosing a suitable $\alpha\in(0,1)$ (see \cite[Lemma 3.2]{scrimalli}). In fact, $\mathcal{K}(\phi)=\nu(\phi)+\mathcal{M}_\circ$, where $\nu(\phi)=(0,\cdots,0,\alpha p_i^0,\cdots,\alpha p_i^0)_{i=1}^{10} $ is $\alpha$-Lipschitz map and $\mathcal{M}_\circ=\{(a_i,\cdots,a_i,b_i,\cdots,b_i,0,\cdots,0)_{i=1}^{10}|\, a_i=A_i,b_i\in [\underbar b_i,\bar b_i]\}$.
Note that for performing this algorithm, we have to consider $p_i^0$ and $c_i$ are no more variables. By taking the demand $D=1$, for $y=(A_i,\cdots,A_i,b_i,\cdots,b_i,\alpha p_i^0,\cdots,\alpha p_i^0)_{i=1}^{10}$ and $z=(A_i,\cdots,A_i,v_i,\cdots,v_i,\alpha p_i^0,\cdots,\alpha p_i^0)_{i=1}^{10}$, the bifunction $f$ can be defined as
\begin{equation}\label{bifunction_game}
	f(y,z)=\sum_{i=1}^{10}\langle -\nabla_{b_i}\theta_i(b),v_i-b_i\rangle,
\end{equation}
 where $\theta_i(b)=(b_i-B_i) q_i(b)+\alpha p_i^0$.
 
 In the view of Corollary \ref{Coro_GNEP_Str_mono}, we apply Algorithm \ref{Algo3} to find the projected solution of EMM for $10$ producers with random values of $A_{i}\in [5,10]$  and $B_{i}\in [-5,5],$ $i=1,2,\ldots,10$ taking different initial points $y_{0}$. Furthermore, the bounds of $b_{i}$ are taken as follows:  $b_1 \in [-1,1]$, $b_2 \in [0,2]$, $b_3 \in [1,3]$, $b_4 \in [-2,-1]$, $b_5 \in [-1,0]$, $b_6 \in [2,3]$, $b_7 \in [1,2]$, $b_8 \in [0.5,1]$, $b_9 \in [-0.5,0]$ and $b_{10} \in [0,50]$ and the values of $p^0_{i}$ are chosen randomly from the interval $[5,10]$ for all $i=1,2,\cdots,10$ and for each initial points. Finally, the bounds of $p_i^1$ for all $i=1,2,\cdots,10$ are taken as follows: $p_1^1 \in [15,40]$, $p_2^1 \in [20,60]$, $p_3^1 \in [17,48]$, $p_4^1 \in [21,51]$, $p_5^1 \in [80,100]$, $p_6^1 \in [30,60]$, $p_7^1 \in [25,53]$, $p_8^1 \in [19,32]$, $p_9^1 \in [25,60]$ and $p_{10}^1 \in [60,90]$. The results obtained by taking different initial points are shown in Table \ref{EMM_table2} with accuracy level $\varepsilon = 10^{-4}$. In Table \ref{EMM_table2}, the symbol \textbf{e} denotes the vector of all $1$'s of dimension 10 and $\mu$ denotes the strong monotonicity parameter for the bifunction $f$ defined in \eqref{bifunction_game}. Since all the conditions in the Corollary \ref{Coro_GNEP_Str_mono} are satisfied, we can ensure that the second last column in Table \ref{EMM_table2} represents the projected solution for the EMM.

% \begin{table}[]
%     \centering
%     \begin{tabular}{c|c}
%       Bounds on $p_i^1$   & Bounds on $b_i$  \\
%          & 
%     \end{tabular}
%     \caption{Caption}
%     \label{tab:my_label}
% \end{table}

% A constraint map $\mathcal{K}$ satisfying (\ref{map_K}) is also known as moving set constraint map. This type of map fulfills hypothesis \ref{hy_3.1ii} of Theorem \ref{Thm_strong_mono}. In the electricity market model, if we have $K_i:I(C)\rightrightarrows \mathbb{R}^{3k}$ as,
%\begin{equation}
%	K_i(\phi)=\{y_i|y_i(q)=a_i q_i^2+b_i q_i+c_i, a_i=A_i,b_i\in [\underbar b_i,\bar b_i],c_i=\alpha p_i^0\}
%\end{equation}
%for some $\alpha\in(0,1)$. Then, $K_i$ is the moving set constraint map satisfying hypothesis (iii). In fact, $K_i(\phi)=\nu(\phi)+A_\circ$, where $\nu(\phi)=(0,\cdots,0,\alpha p_i^0,\cdots,\alpha p_i^0)$ and $A_\circ=\{(a_i,\cdots,a_i,b_i,\cdots,b_i,0,\cdots,0)|\, a_i=A_i,b_i\in [\underbar b_i,\bar b_i]\}$.

\begin{table}[!h]
\centering
\caption{Results obtained by Algorithm \ref{Algo3} to electricity market model for $10$ producers}
\label{EMM_table2}
\renewcommand{\arraystretch}{0.5}
%{\scriptsize
%\resizebox{\columnwidth}{!}{%
\begin{tabular}{p{1cm}|p{0.3cm}|p{1.2cm}|p{0.7cm}|p{0.6cm}|p{2.9cm}|p{3.2cm}|p{1.5cm}}
\hline
 Initial point $y_{0}=(b_i^0)_{i=1}^{10}$ & iter &time (in sec) &$\mu$&$L$&Obtained solution $\bar{y}=(b_{i})_{i=1}^{10}$&$(p_{i}^{1})_{i=1}^{10}$ & Conclusion \\ \hline
  \textbf{e} & $4$   & $1668.61$ & $0.020$& $0.17$ &(0.99,0,1,$-$1.30,0,2, 1,0.94,0,$-$2)& (15.24,21.25,17.26,80, 30,25,25.12,25,60)& (Projected solution)   \\  \hline
   $2$ \textbf{e} &  $6$  & $2174.29$& $0.021$& $0.18$  &(0.08,0,1,$-$1,0,2,1, 0.99,0,$-$2)& (15,20,17,21,80,30,25, 20.01,25,60)  & (Projected solution)    \\  \hline
    $ 3$ \textbf{e} &  $4$  & $1393.03$& $0.023$& $0.18$  &(0.34,1.99,1,$-$1,0, 2.69,1.99,0.99,0,$-$2)& (15.09,25.19,21,80,30, 26.86,19,25,60)  & (Projected solution)    \\  \hline
      $ -\frac{1}{2}$ \textbf{e} &  $5$  & $2051.61$& $0.018$& $0.16$  &($-$0.56,0.71,2.99,$-$1, 0,2,1,0.99,0,$-$2)& (15,20,25.77,21.01,80, 30,25.91,25.05,25,60)& (Projected solution)    \\ \hline
\end{tabular}
\end{table}

%   \begin{table}[h!]
%       \caption{Test Table}
%   \begin{adjustbox}{max width=\textwidth}
%   \begin{tabular}{*{14}{|c}|}%%{|c|c|c|c|c|c|c|c|c|c|c|c|c|c|}
%   \hline
%   One & Two &Three & Four & Five & Six & Seven & Eight & Nine & Ten & Eleven &
%   Twelve & Thirteen & Fourteen\\
%   \hline
%   \hline
%   $1.111$ & $2.222$ & $3.333$ & $4.444$ & $5.555$ & $6.666$ & $7.777$ &
%   $8.888$ & $9.999$ & $0.000$ & $1.111$ & $2.222$ & $3.333$ & $4.444$\\
%   \hline
% \end{tabular}
% \end{adjustbox}
% \end{table}

%$$\widetilde{\varPi}_F(0)$$

\section{Conclusion}
In this paper, two iterative algorithms have been developed to find the projected solution of QEP. The first algorithm is suitable for solving the strongly monotone QEP, while the other is designed for the monotone QEP. We proved the convergence of the sequences generated by both algorithms to the projected solution of the QEP. We constructed some examples based on the assumptions for which the convergence of the generated sequence of proposed algorithms was proved, and the performance of both algorithms was numerically tested. In addition, we solved the electricity market model with two and ten producers with the help of the proposed algorithms.

\end{document}